\documentclass[11pt]{article}
\usepackage{mathrsfs}
\usepackage{}
\usepackage{amsfonts}
\usepackage{amssymb}
\usepackage{graphicx}
\usepackage{amsmath}

\def\sqr#1#2{{\vcenter{\vbox{\hrule height.#2pt
              \hbox{\vrule width.#2pt height#1pt \kern#1pt \vrule
width.#2pt}
              \hrule height.#2pt}}}}
\def\signed #1{{\unskip\nobreak\hfil\penalty50
              \hskip2em\hbox{}\nobreak\hfil#1
              \parfillskip=0pt \finalhyphendemerits=0 \par}}
\def\endpf{\signed {$\sqr69$}}
\def\3n{\negthinspace \negthinspace \negthinspace }
\def\2n{\negthinspace \negthinspace }
\def\1n{\negthinspace }

\def\ms{\medskip}

\def\({\Big (}
\def\){\Big )}
\def\[{\Big[}
\def\]{\Big]}

\def\be{\begin{equation}}
\def\bel{\begin{equation}\label}
\def\ee{\end{equation}}
\def\bea{\begin{eqnarray}}
\def\eea{\end{eqnarray}}
\def\bt{\begin{theorem}}
\def\et{\end{theorem}}
\def\bc{\begin{corollary}}
\def\ec{\end{corollary}}
\def\bl{\begin{lemma}}
\def\el{\end{lemma}}
\def\bp{\begin{proposition}}
\def\ep{\end{proposition}}
\def\br{\begin{remark}}
\def\er{\end{remark}}
\def\ba{\begin{array}}
\def\ea{\end{array}}
\def\bd{\begin{definition}}
\def\ed{\end{definition}}

\newtheorem{lemma}{Lemma}[section]
\newtheorem{remark}{Remark}[section]

\newtheorem{theorem}{Theorem}[section]
\newtheorem{corollary}{Corollary}[section]

\newtheorem{definition}{Definition}[section]
\newtheorem{proposition}{Proposition}[section]

\makeatletter
   
   \@addtoreset{equation}{section}
\makeatother
\allowdisplaybreaks
\begin{document}

\title{\bf  Local controllability of a free-boundary problem for a class of one-dimensional degenerate parabolic equations
 }

\author{Lingyang Liu\thanks{School of Mathematical Sciences, South China Normal
University, Guangzhou, 510631, China. E-mail address:
liuly@scnu.edu.cn   \ms } }

\date{}

\maketitle

\begin{abstract}
This paper is devoted to a study of the controllability of a free-boundary problem for a class of one-dimensional degenerate parabolic equations with distributed controls, locally supported in space. We prove that for any  $T>0$, if the initial state is sufficiently small, there exists a control that drives the state exactly to rest at
time $t = T$. The proof is based on Schauder's fixed point theorem, combined with appropriate estimates for solutions to degenerate parabolic equations and for control functions.
\end{abstract}

\noindent{\bf Key Words.} Controllability, Carleman estimate, free-boundary problem,
1D degenerate parabolic equation

\section{Introduction and main result}

In real-world applications, many physical phenomena are modeled by degenerate parabolic equations, such as population genetics, fluid dynamics, and vision (see \cite{CMV} for a more detailed introduction). Due to these wide-ranging applications, control problems for such equations have received significant attention. Different from nondegenerate parabolic equations, a key feature of degenerate problems is that the diffusion coefficient in the principal operator fails to maintain uniform ellipticity. Consequently, new tools have been developed to study related controllability and observability issues. For instance, novel Carleman-type estimates with adapted weight functions, distinct from those for nondegenerate parabolic equations, have been employed to derive observability inequalities for the corresponding adjoint problems (see, e.g., \cite{ACF,4}).

To date, control issues concerning different aspects of degenerate parabolic equations in cylindrical domains have been extensively studied. The investigation mainly focuses on the one-dimensional case. Gueye \cite{Gue} investigated the boundary null controllability of the one-dimensional degenerate heat equation by the transmutation method, when the control acts through the part of the boundary where degeneracy occurs.  Cannarsa et al. \cite{5}  discussed the null controllability of the one-dimensional degenerate heat equation using spectral analysis and the moment method, when the degeneracy occurs inside the domain.  In addition, Ara\'{u}jo et al. \cite{ADV} recovered the boundary null controllability for the degenerate heat equation by analyzing the asymptotic behaviour of an eligible family of state-control pairs and solving the corresponding singularly perturbed internal null controllability problems. Regarding other aspects, Cannarsa et al. \cite{CMV1} studied the null controllability cost for a one-dimensional weakly degenerate parabolic equation governed by a boundary control acting at the degeneracy point  by means of the moment method. Later, they established a result for strongly degenerate parabolic equations controlled either by a boundary control or by a locally distributed control in \cite{CMV2}. Besides, Benoit et al. \cite{BLR} investigated the null controllability of a class of linear one-dimensional strongly degenerate parabolic equations with measurable coefficients by applying the flatness approach. We also mention that Wang and Du \cite{WD1} proved the null controllability of a class of semilinear one-dimensional weakly degenerate parabolic equations involving convection terms. Their approach relies on Carleman estimates for the linearized system and a fixed-point argument.

As we know, controllability of degenerate parabolic equations is dependent on the degree of degeneracy. If the degeneracy is too strong, systems may not be null controllable. In such cases, Wang investigated the approximate controllability of a class of linear and semilinear degenerate systems in \cite{Wang1} and \cite{Wang}, respectively. In \cite{WD}, Wang and Du also established the approximate controllability result for a class of semilinear degenerate systems where the convection terms are dominated by the diffusion terms.

In practice, many dynamical processes evolve in domains with moving boundaries. A classic example is the changing interface of an ice-water mixture under rising temperatures. Over the last few decades, researchers have studied the stability of degenerate parabolic equations in non-cylindrical domains. For example, M. Bertsch et al. \cite{BDF} analyzed the behavior of solutions near the vertex $(0,0)$ for a degenerate parabolic equation in a special class of non-cylindrical domains. In \cite{GLL}, Gao et al. investigated the asymptotic behavior of solutions to the degenerate heat equation in non-cylindrical domains with linear moving boundaries. However, to the best of our knowledge,
controllability problems for degenerate parabolic equations with moving boundaries (or a free-boundary) have not been addressed.

Let $T>0$, $L_0>0$ and $0<a<b<L_\ast<L_0$ be given. For any function $L\in C^1([0,T])$ with $L(t)>0$, set $Q_L:=\{(x,t):x\in(0,L(t)),t\in(0,T)\}$. Write $\omega=(a,b)$ and denote by $\chi_\omega$ the characteristic function of $\omega$. For any $0 \leq \alpha < 2$, we consider a free-boundary problem for degenerate parabolic systems of the form
\begin{equation}\label{1.1}
\left\{\begin{array}{ll}
y_{t}-(x^\alpha y_x)_x=1_\omega v,& (x,t)\in Q_L,\\[2mm]
\left\{\begin{array}{ll}
y(0,t)=0, & \text{if}\quad \alpha\in[0,1),\\[2mm]
\lim\limits_{x\rightarrow0^+}x^\alpha y_x(x,t)=0, &\text{if}\quad \alpha\in[1,2),\\[2mm]
\end{array}\right.&\\[2mm]
y(L(t),t)=0, &t\in(0,T),\\[2mm]
y(x,0)=y_0(x),& x\in(0,L_0),\\[2mm]
L(0)=L_0,
\end{array}\right.
\end{equation}
with the additional boundary condition
\begin{equation}\label{1.2}
L'(t)=-L^\alpha(t)y_x(L(t),t),\quad t\in(0,T).
\end{equation}
Here, $y$ is the state variable and $v$ is a control; $v$ acts on the system at any time via the nonempty open set $\omega$.

As observed, \eqref{1.1} degenerates at $x=0$ for $\alpha>0$. More precisely, the degeneracy is weak if $0<\alpha<1$, while strong if $\alpha\geq1$. Note that in \eqref{1.1}, a portion of the domain boundary is not predetermined but evolves as part of the solution (referred to as a free boundary), which differs fundamentally from non-cylindrical domain problems where boundary evolution is prescribed. Therefore, to identify the state $(L,y)$, additional information is needed (termed a free boundary condition) and this is provided by \eqref{1.2}.
In many areas, this condition is completely natural. For example, consider dissipative systems derived from conservation laws:
\begin{equation}\label{1.3}
y_t-\text{div}(k(x,t)\nabla y)=0,
\end{equation}
where $y$ denotes the temperature and $k\nabla y(x,t)$ represents the heat flux. If $k = 1$, obviously, \eqref{1.3} becomes the classical heat equation and a free boundary condition is given by
\begin{equation*}
L'(t)=-y_x(L(t),t),\quad t\in(0,T).
\end{equation*}
Whereas, if the heat
conduction coefficient $k$ depends on the degeneracy parameter in a manner as $x^\alpha$, then \eqref{1.2} is proposed. In particular, when $\alpha=0$, \eqref{1.2} reduces to the above boundary condition.

It is well known that free-boundary problems for parabolic equations serve as mathematical models for a variety of phenomena in nature, such as phase transitions in materials science \cite{Fri}, tumor growth in biomedical modeling \cite{Fri1}, and gas flow through porous media \cite{Vaz}, among others. In practical applications, \eqref{1.1} may describe the dynamical behavior of a diffusion system where, at one end of the material ($ x=0 $), the diffusion coefficient approaches zero, reflecting the cessation of heat or mass transfer at that boundary (e.g., due to perfect insulation or material depletion), while the interface position $L(t)$ evolves over time, governed by energy or mass conservation laws.

Let us mention some previous works on the controllability of free-boundary problems for non-degenerate parabolic equations. In \cite{FLde} and \cite{Fde}, Fern\'{a}ndez-Cara et al. investigated the local null controllability of free boundary problems for the linear and semilinear 1D heat equations, respectively. Wang et al. \cite{WLL} considered a free boundary problem for a class of quasi-linear 1D parabolic equations and established a local null controllability result. V. Costa et al. \cite{CLLP} analyzed the local null controllability of a free-boundary problem for the 1D heat equation with local and nonlocal nonlinearities. Concerning the controllability to the trajectories of the one-phase Stefan problem, J. A. B\'{a}rcena-Petisco et al. \cite{BFS} proved that the state is locally controllable to smooth trajectories with nonnegative controls. B. Colle et al. \cite{CLT} studied the controllability of the Stefan problem with two phases and provided numerical illustrations of their results. We also note the work by B. Geshkovski and E. Zuazua \cite{GZ}, who addressed the local controllability of a one-dimensional free boundary problem for a fluid governed by the viscous Burgers equation.

In order to study the well-posedness of the degenerate parabolic problem, we introduce some weighted Sobolev spaces. Let $L\in C^1([0,T])$ be given. For $0\leq\alpha<1$, define the Hilbert space $H_\alpha^1(0,L(t))$ as follows:
\begin{eqnarray*}
&&H_\alpha^1(0,L(t)):=\big{\{}y\in L^2(0, L(t))\ \big{|}\ y\ \text{is absolutely\
continuous\ in}\ [0,L(t)],\\
&&\quad\quad\quad\quad\quad\quad\quad\quad\quad\quad\quad\quad\quad\quad x^\frac{\alpha}{2} y_x\in L^2(0, L(t))\
\text{and}\ y(0)=y(L(t))=0\big\}.
\end{eqnarray*}
In addition, set
\begin{eqnarray*}
D(0,L(t)):=\big{\{}y\in H_\alpha^1(0,L(t))\ \big{|}\ x^\alpha y_x\in H^1(0,L(t))\big\}.
\end{eqnarray*}
Notice that if $y \in D(0,L(t))$ (or even $y \in H^1_\alpha(0,L(t))$, then $y$ satisfies the Dirichlet boundary
conditions $y(0) = y(L(t)) = 0$.

For $1\leq\alpha<2$, $H_\alpha^1(0,L(t))$ is defined as
\begin{eqnarray*}
&&H_\alpha^1(0, L(t)):=\big{\{}y\in L^2(0, L(t))\ \big{|}\ y\ \text{is locally\ absolutely\
continuous\ in}\ (0, L(t)],\\
&&\quad\quad\quad\quad\quad\quad\quad\quad\quad\quad\quad\quad\quad\quad x^\frac{\alpha}{2} y_x\in L^2(0, L(t))\ \text{and}\ y(L(t))=0\big\}.
\end{eqnarray*}
Then
\begin{eqnarray*}
&&D(0,L(t)):=\{y \in H^1_\alpha(0,L(t))\ \big{|}\ x^\alpha y_x \in H^1(0,L(t))\}\\[2mm]
&&\quad\quad\quad\quad\quad =\big{\{}y\in L^2(0, L(t))\ \big{|}\ y\ \text{is locally\ absolutely\
continuous\ in}\ (0, L(t)],\\[2mm]
&&\quad\quad\quad\quad\quad\quad\quad\quad\quad x^{\alpha} y\in H^1_0(0, L(t)),\ x^{\alpha} y_x\in H^1(0,L(t))\ \text{and}\ \lim\limits_{x\rightarrow0^+}x^\alpha y_x(x)=0 \big\}.
\end{eqnarray*}
Notice that if $y \in D(0,L(t))$, then $y$ satisfies the Neumann boundary condition $\lim\limits_{x\rightarrow0^+}x^\alpha y_x(x)=0$ and the Dirichlet boundary condition $y(L(t)) = 0$ at $x = L(t)$.

It is well known that $H_\alpha^1(0,L(t))$ is a Banach space endowed with the norm
\begin{eqnarray*}
\|y\|_{H_\alpha^1(0,L(t))} = \|y\|_{L^2(0,L(t))} + \|x^{\alpha/2}y_x\|_{L^2(0,L(t))}\quad\ \forall y\in H_\alpha^1(0,L(t)).
\end{eqnarray*}

The goal of this paper is to study the null controllability of \eqref{1.1}--\eqref{1.2}. To this end, we formulate the following controllability concepts:
\begin{definition}
Problem \eqref{1.1}--\eqref{1.2} is said to be null-controllable at time $T$, if for any initial value $y_0\in H^1_\alpha(0,L_0)$, there exist a control $v\in L^2(\omega\times(0,T))$,
a function $L\in C^1([0,T])$ and an associated solution
$y$ satisfying \eqref{1.1}--\eqref{1.2} and
\begin{equation}\label{d1}
y(x,T)=0,\quad x\in(0,L(T)).
\end{equation}
\end{definition}

\begin{definition}
Problem \eqref{1.1}--\eqref{1.2} is said to be approximately null-controllable at time $T$, if for any initial value $y_0\in H^1_\alpha(0,L_0)$ and any $\beta>0$, there exist a control $v\in L^2(\omega\times(0,T))$,
a function $L\in C^1([0,T])$ and an associated solution
$y$ satisfying \eqref{1.1}--\eqref{1.2} and
\begin{equation}\label{d2}
\|y(\cdot,T\|_{L^2(0,L(T))}\leq\beta.
\end{equation}
\end{definition}
\begin{remark}
Approximate controllability combined with uniform estimates on the approximate controls, as $\beta\rightarrow0$, may lead to null controllability.
\end{remark}

\begin{definition}\label{df}
Problem \eqref{1.1}--\eqref{1.2} is said to be locally null-controllable at time $T$, if there exists $\varepsilon> 0$ such that for any given initial value $y_0 \in H^1_\alpha(0,L_0)$ satisfying $\|y_0\|_{H^1_\alpha(0,L_0)}\leq \varepsilon$, there exist a control $v\in L^2(\omega\times(0,T))$,
a function $L\in C^1([0,T])$ and an associated solution
$y$ satisfying \eqref{1.1}--\eqref{1.2} and \eqref{d1}.
\end{definition}

The main result is stated as follows.
\begin{theorem}\label{t1.1}
Suppose that $0\leq\alpha<2$. For any given $T>0$, problem \eqref{1.1}--\eqref{1.2} is locally null-controllable at time $T$ in the sense of Definition \ref{df}.
\end{theorem}

\begin{remark}\label{re1}
In this paper, we discuss the null controllability of \eqref{1.1}--\eqref{1.2} under the parameter condition $0\leq\alpha<2$. This is because when $\alpha\geq2$, equation \eqref{1.1} may not be null controllable. Indeed, \cite{CMV0} provides an explicit example demonstrating the failure of null controllability for degenerate parabolic equations in a cylindrical domain when $\alpha\geq2$.
\end{remark}

\begin{remark}
When $\alpha=0$, equation \eqref{1.1} reduces to the classical heat equation and the assumption on the initial condition becomes $y_0\in H^1_0(0,L_0)$ with $\|y_0\|_{H^1_0(0,L_0)}\leq \varepsilon$. In this case, our results align with earlier findings for the classical heat equation (see \cite{FLde}).
\end{remark}

\begin{remark}
It would be interesting to establish a global controllability result for
\eqref{1.1}--\eqref{1.2}. However, by means of the fixed point
technique (detailed in Section \ref{sec3}), this seems very difficult and it remains to be an unsolved problem (it is not clear how the
fixed point mapping $\Lambda_\beta$ may be pulled back from $\mathscr M$ into itself when the initial value $y_0$ is large).
On the other hand, we note that even for the heat equation, the existence of a global-in-time solution fails for general initial data $y_0$ in free boundary problems (see \cite{FP}).
\end{remark}

\begin{remark}
In multidimensional settings, a challenging issue in free boundary problems lies in the study of regularity and certain geometric properties of free boundaries (see, e.g., \cite{DF}). Controllability of free boundary problems for degenerate parabolic equations in higher dimensions remains open.
\end{remark}

The proof strategy for Theorem \ref{t1.1} is outlined below:

1. We first prove that \eqref{1.1}--\eqref{1.2} is uniformly approximately null-controllable, i.e., for any $\beta>0$, there exists a
triplet $(L_\beta, y_\beta,v_\beta)$ is uniformly bounded in an appropriate space and satisfy \eqref{1.1}--\eqref{1.2} and \eqref{d2}. For this purpose, we introduce a fixed-point reformulation based on suitable linearized problems and verify that Schauder's theorem can be applied if the initial data $y_0$ is sufficiently small.

Note that although the degeneracy occurs at the left fixed endpoint, its effects propagate across the entire domain. In order to analyze the H\"{o}lder regularity of the controlled solution near the free boundary, we need to prove the existence, uniqueness and regularity of solutions to the
degenerate parabolic systems in non-cylindrical domains and derive some necessary compactness estimates for solutions. Moreover, a controllability result for the linearized system, along with cost estimates for the control functions, is required.

2. Then, by taking the limit as $\beta\rightarrow0$, we see that such triplets (at least for a subsequence) converge to a solution of \eqref{1.1}--\eqref{1.2} and \eqref{d1}.

From a technical perspective, problem \eqref{1.1}--\eqref{1.2} presents greater challenges than those involving only boundary degeneracy or a free boundary, for reasons we now elucidate:

To prove the well-posedness of \eqref{1.1}, we need to address both the degeneracy and time-dependent coefficients in the parabolic operator by employing a diffeomorphism that transforms \eqref{1.1} into an equivalent cylindrical parabolic problem, which requires more intricate calculations and delicate estimates. On the other hand, the controllability of system \eqref{1.1} is not only determined by the parameter $\alpha$ but also by the moving boundary. To solve the controllability problem, we shall derive a new Carleman estimate to establish an observability inequality for the corresponding dual system. This requires a careful analysis of weighted integrals (in the domain and on the boundary) of solutions to degenerate parabolic problems. Finally, by selecting appropriate weighted functions and applying Hardy-type inequalities, we obtain the desired Carleman estimate.

The rest of the paper is organized as follows. In Section 2, we demonstrate the well-posedness and controllability of degenerate parabolic equations in non-cylindrical domains. Section 3 is devoted to proving Theorem \ref{t1.1}. In Section 4, we derive a global Carleman estimate stated in Theorem \ref{t2.2}. Finally, Section 5 serves as an appendix providing the proofs of two kinds of useful inequalities for non-cylindrical domains.

\section{Some results for degenerate parabolic equations in non-cylindrical domains}\label{s2}

In this section, we discuss the well-posedness and controllability of the controlled linear system  \eqref{1.1} (excluding \eqref{1.2}), under the assumption that $L\in C^1([0,T])$ satisfies $L(0)=L_0$ and $L_\ast\leq L(t)\leq B$ for all $t\in[0,T]$. The results on the existence and uniqueness of solutions, along with some compactness estimates, will play a fundamental role in solving free boundary problems by a fixed point method.

For simplicity, we set $N_L:=\|L'\|_\infty$ and let $\omega_0$ be a non-empty open subset of $\omega=(a,b)$ such that $\overline{\omega_0}\subseteq\omega$.

\subsection{Well-posedness of the degenerate parabolic problem and some compact estimates}
In general, we consider the initial-boundary value problem
\begin{equation}\label{2.1}
\left\{\begin{array}{ll}
y_t-(x^\alpha y_x)_x=f, & (x,t)\in Q_L,\\[2mm]
\left\{\begin{array}{ll}
y(0,t)=0, & \text{if}\quad \alpha\in[0,1),\\[2mm]
\lim\limits_{x\rightarrow0^+}x^\alpha y_x(x,t)=0, &\text{if}\quad \alpha\in[1,2),\\[2mm]
\end{array}\right.&\\[2mm]
y(L(t),t)=0, &t\in(0,T),\\[2mm]
y(x,0)=y_0(x), &x\in (0,L_0),
\end{array}\right.
\end{equation}
where $y_0\in L^2(0,L_0)$ and $f\in L^2(Q_L)$.

First, we define  the weak solution to problem \eqref{2.1} as follows.
\begin{definition}
A function $y$ is called a weak solution of problem \eqref{2.1}, if $y\in C([0,T];L^2(0,L(t)))\cap L^2(0,T;H_\alpha^1(0,L(t)))$ and for any function $\varphi\in L^\infty(0,T;L^2(0,L(t)))\cap L^2(0,T;H_\alpha^1(0,L(t)))$ with $\varphi_t\in L^2(Q_L)$ and $\varphi(\cdot,T)|_{(0,L(T))}=0$, the following integral equality holds:
\begin{equation*}\label{d2.1}
\iint_{Q L}-y\varphi_t+x^\alpha y_x \varphi_x dx dt-\int_{0}^{L_0}y_{0}(x)\varphi(x,0)dx=\iint_{Q_L}f\varphi dx dt.
\end{equation*}
\end{definition}

Then we have the following well-posedness result for \eqref{2.1}:
\begin{theorem}\label{t2.1}
$(i)$ Let $f$ be given in $L^2(Q_L)$. For any $y_0\in L^2(0,L_0)$, problem \eqref{2.1} admits a unique weak solution
$
y\in C([0,T];L^2(0,L(t)))\cap L^2(0,T;H_\alpha^1(0,L(t))).
$
Moreover, the solution $y$ satisfies the estimate
\begin{equation*}
\|y\|^2_{L^\infty(0,T;L^2(0,L(t)))}+\|x^\alpha y_x^2\|_{L^1(Q_L)}\leq e^{T(\|L^{-1}L'\|_\infty+1)}\big( \|f\|^2_{L^2(Q_L)}+\|y_0\|^2_{L^2(0,L_0)}\big).
\end{equation*}

\item$(ii)$ If $y_0\in H^1_\alpha(0,L_0)$ additionally, then $y_t\in L^2(Q_L)$, $x^\alpha y_x^2\in L^\infty(0,T;L^1(0,L(t)))$ and $y\in L^2(0,T;D(0,L(t)))$. It holds that
\begin{eqnarray*}\label{2.7}
\|y_t\|^2_{L^2(Q_L)}+\|x^\alpha y_x^2\|_{L^\infty(0,T;L^1(0,L(t)))}\leq\!\!\!\!\!\!\!\!\!\!&& 4e^{2T(\|L^{-1}L'\|^2_\infty\|L^{2-\alpha}\|_\infty+(2-\alpha)L_0\|L^{-1}L'\|_\infty)}\\[2mm]
&&\quad \times\big( 2\|f\|_{L^2(Q_L)}+\|x^\alpha y_{\scriptscriptstyle0,x}^2\|_{L^1(0,L_0)}\big)
\end{eqnarray*}
and
\begin{eqnarray*}
&&\big\|\big(x^\alpha y_x\big)_x\big\|^2_{L^2(Q_L)}\leq 16\big(\|L^{-1}L'\|^2_\infty\|L^{2-\alpha}\|_\infty+1\big)\notag\\[2mm]
&&\quad\times e^{2T(\|L^{-1}L'\|^2_\infty\|L^{2-\alpha}\|_\infty+(2-\alpha)L_0\|L^{-1}L'\|_\infty)}
\big(3\|f\|^2_{L^2(Q_L)}+\|x^\alpha |y_{\scriptscriptstyle0,x}|^2\|_{L^1(0,L_0)}\big).
\end{eqnarray*}

\item$(iii)$ If $y_0\in L^\infty(0,L_0)$ and $f\in L^\infty(Q_L)$ additionally, then $y\in L^\infty(Q_L)$ and
\begin{equation*}
\|y\|_{L^\infty(Q_L)}\leq T\|f\|_{L^\infty(Q_L)}+\|y_0\|_{L^\infty(0,L_0)}.
\end{equation*}
\end{theorem}

To prove Theorem \ref{t2.1}, we first need to reformulate the problem. Specifically, (\ref{2.1}) can be transformed into an equivalent cylindrical system via the following change of variables:
\begin{equation}\label{2.2}
\Phi_L: Q_L\rightarrow Q,\quad\quad \Phi_L(x,t)=(\varsigma,t)=\(\text{\footnotesize$\frac{  L_{\scriptscriptstyle0}x}{ L(t)}$}, \text{\small$t$}\),
\end{equation}
where $Q:=(0,L_0)\times(0,T)$. Direct calculation shows that the Jacobian determinant associated with this change of variables satisfies
\begin{equation*}
J(\varsigma,t)=
\left|\begin{array}{ll}
x_\varsigma & x_t\\
0 &1
\end{array}\right|=\frac{\scriptstyle L(t)}{\scriptstyle L_0}>\frac{\scriptstyle L_\ast}{\scriptstyle L_0}>\text{\footnotesize $0$},\quad\ (\varsigma,t)\in \overline{Q}.
\end{equation*}
Define $w(\varsigma,t):=y(x,t)$. Then $w$ is well-defined in the rectangle $\overline{Q}$ and by \eqref{2.2} it is easy to verify that
$$y_t=w_t-\frac{\scriptstyle L'(t)}{\scriptstyle L(t)}\varsigma w_\varsigma,\quad\
y_x=\frac{\scriptstyle L_0}{\scriptstyle L(t)}w_\varsigma,\quad\
(x^\alpha y_x)_x=\text{$\(\frac{\scriptstyle L(t)}{\scriptstyle L_0}\)$}^{\scriptscriptstyle\alpha-2}(\varsigma^\alpha w_\varsigma)_\varsigma.
$$
Thus, $w$ solves the problem
\begin{equation}\label{2.3}
\left\{\begin{array}{ll}
w_{t}-p(t)(\varsigma^\alpha w_\varsigma)_\varsigma- q(t)\varsigma w_\varsigma=h, & (\varsigma,t)\in Q,\\[2mm]
\left\{\begin{array}{ll}
w(\text{\small$0$},t)=0, & \text{if}\quad \alpha\in[0,1),\\[2mm]
\lim\limits_{\varsigma\rightarrow0^+}(\varsigma^\alpha w_\varsigma)(\varsigma,t)=0, &\text{if}\quad \alpha\in[1,2),\\[2mm]
\end{array}\right.&\\[2mm]
w(\text{\small $L_0$},t)=0, &t\in(0,T),\\[2mm]
w(\varsigma,\text{\small$0$})=w_{\scriptscriptstyle0}(\varsigma), &\varsigma\in (0,L_0),
\end{array}\right.
\end{equation}
where $\text {\small $p(t)$}=\big(\frac{\scriptscriptstyle L(t)}{\scriptscriptstyle L_0}\big)^{\scriptscriptstyle \alpha-2}$ and $\text{\small $q(t)$}=\frac{\scriptscriptstyle L'(t)}{\scriptscriptstyle L(t)}$.

To establish the well-posedness result in Theorem \ref{t2.1}, we need only to prove that problem \eqref{2.3} is well-posed and that certain estimates hold, since $\Phi_L$ is reversible. For this purpose, we define the weak solution to problem \eqref{2.3}:
\begin{definition}
A function $w$ is called a weak solution of the problem \eqref{2.3}, if $w\in C([0,T];L^2(0,L_0))\cap L^2(0,T;H_\alpha^1(0,L_0))$ and for any function $\psi\in L^\infty(0,T;L^2(0,L_0))\cap L^2(0,T;H_\alpha^1(0,L_0))$ with $\psi_t\in L^2(Q)$ and $\psi(\cdot,T)|_{(0,L_0)}=0$, the following integral equality holds:
\begin{equation*}\label{d2.1}
\iint_{Q}-w\psi_t+p(t)\varsigma^\alpha w_\varsigma \psi_\varsigma-q(t)\varsigma w_\varsigma \psi d\varsigma dt-\int_{0}^{L_0}w_{0}(\varsigma)\psi(\varsigma,0)d\varsigma=\iint_{Q}h\psi d\varsigma dt.
\end{equation*}
\end{definition}
\begin{remark}
Assume that $w\in C([0,T];L^2(0,L_0))\cap L^2(0,T;H_\alpha^1(0,L_0))$ with $w_t\in L^2(Q)$. Then $w$ is a weak solution of the problem \eqref{2.3}, if and only if the integral equality
\begin{equation*}
\iint_{Q}w_t\psi+p(t)\varsigma^\alpha w_\varsigma \psi_\varsigma-q(t)\varsigma w_\varsigma \psi d\varsigma dt=\iint_{Q}h\psi d\varsigma dt
\end{equation*}
holds for any function $\psi\in L^\infty(0,T;L^2(0,L_0))\cap L^2(0,T;H_\alpha^1(0,L_0))$, and $w(\cdot,\text{\small$0$})=w_{\scriptscriptstyle0}$ in the trace sense.

\end{remark}
Problem \eqref{2.3} is well-posed. We have
\begin{proposition}\label{p2.1}
$(i)$ For any $h\in L^2(Q)$ and $w_0\in L^2(0,L_0)$, problem \eqref{2.3} admits a unique weak solution $w\in C([0,T];L^2(0,L_0))\cap L^2(0,T;H_\alpha^1(0,L_0))$. Furthermore, the solution $w$ satisfies
\begin{equation*}\label{2.6}
\|w\|^2_{L^\infty(0,T;L^2(0,L_0))}+\|p(t)\varsigma^\alpha w_\varsigma^2\|_{L^1(Q)}\leq  e^{T(\|q\|_{\infty}+1)}\big(\|h\|^2_{L^2(Q)}+\|w_0\|^2_{L^2(0,L_0)}\big).
\end{equation*}

\item$(ii)$ If $w_0\in H^1_\alpha(0,L_0)$ additionally, then $w_t\in L^2(Q)$, $p(t)\varsigma^\alpha w_\varsigma^2\in L^\infty(0,T;L^1(0,L_0))$ and
$p(t)\big(\varsigma^\alpha w_\varsigma\big)_\varsigma\in L^2(Q)$. Moreover, the following estimates hold:
\begin{eqnarray*}\label{2.7}
\| w_t\|^2_{L^2(Q)}+\|p(t)\varsigma^\alpha w_\varsigma^2\|_{L^\infty(0,T;L^1(0,L_0))}\leq\!\!\!\!\!\!\!\!\!\!&& 4e^{2T({{\scriptstyle L}^{\scriptscriptstyle2-\alpha}_{\scriptscriptstyle0}}\|q\|^2_{\infty}\|\frac{1}{p}\|_{\infty}+\|\frac{p'}{p}\|_{\infty})}\\[2mm]
&& \times\big(2\|h\|^2_{L^2(Q)}+\|\varsigma^\alpha |w_{\scriptscriptstyle0,\varsigma}|^2\|_{L^1(0,L_0)}\big)
\end{eqnarray*}
and
\begin{eqnarray*}
&&\big\|p(t)\big(\varsigma^\alpha w_\varsigma\big)_\varsigma\big\|^2_{L^2(Q)}\leq 16\big({\textstyle L}_{\scriptscriptstyle0}^{\scriptscriptstyle2-\alpha}\|q\|^2_{\infty}\big\|{\textstyle\frac{1}{p}}\big\|_{\infty}+1\big)\notag\\[2mm]
&&\quad\times e^{2T({{\scriptstyle L}^{\scriptscriptstyle2-\alpha}_{\scriptscriptstyle0}}\|q\|^2_{\infty}\|\frac{1}{p}\|_{\infty}+\|\frac{p'}{p}\|_{\infty})}
\big(3\|h\|^2_{L^2(Q)}+\|\varsigma^\alpha |w_{\scriptscriptstyle0,\varsigma}|^2\|_{L^1(0,L_0)}\big).
\end{eqnarray*}

\item$(iii)$ If $w_0\in L^\infty(0,L_0)$ and $h\in L^\infty(Q)$ additionally, then $w\in L^\infty(Q)$ and
\begin{equation*}
\|w\|_{L^\infty(Q)}\leq T\|h\|_{L^\infty(Q)}+\|w_0\|_{L^\infty(Q)}.
\end{equation*}
\end{proposition}

We adopt the parabolic regularization method to establish the existence of solutions, while uniqueness is proved via the Holmgren method. The idea of this proof is inspired by \cite{Wang}. Notably, equation \eqref{2.3} differs from that in \cite{Wang} due to the presence of an additional advection term $q(t)\varsigma w_\varsigma$. Therefore, in the subsequent analysis, we focus on addressing this term and derive some necessary compactness estimates for solutions.

\emph{Proof of Proposition \ref{p2.1}.}
To simplify the notation, we set $a(\varsigma,t) = p(t)\varsigma^\alpha$. For any positive integer $k$, choose $a_k, h_k \in C^\infty(\overline{Q})$, $q_k\in C^\infty([0,T])$ and $w^{(k)}_0 \in C^\infty([0,L_0])$,  satisfying
$$
a(\varsigma,t)+\frac{1}{k}\leq a_k(\varsigma,t)\leq a(\varsigma,t)+\frac{2}{k},\quad \big\|{\textstyle\frac{1}{a_k}\frac{\partial a_k}{\partial t}}\big\|_{L^\infty(Q)}\leq\big\|{\textstyle\frac{1}{a}\frac{\partial a}{\partial t}}\big\|_{L^\infty(Q)},
$$
$$
\|q_k\|_{L^\infty(0,T)}\leq\|q\|_{L^\infty(0,T)},\quad \|h_k \|_ {L^2(Q)} \leq \|h\|_ {L^2(Q)},\quad \|w^{(k)}_0\|_{L^2(0,L_0)}\leq\|w_0\|_{L^2(0,L_0)},
$$
$k=1,2,\cdots,$
and
$$
q_k\rightarrow q\ \text{in}\ L^2(0,T), \quad h_k\rightarrow h\ \text{in}\ L^2(Q), \quad w^{(k)}_0\rightarrow w_0\ \text{in}\ L^2(0,L_0),\quad \text{as}\ k\rightarrow\infty.
$$
Further,
$$
\|a_{\scriptscriptstyle k}(\varsigma,0)|w^{(\scriptscriptstyle k)}_{\scriptscriptstyle0,\varsigma}|^2\|_{L^1(0,L_0)}\leq\|\varsigma^\alpha |w_{\scriptscriptstyle0,\varsigma}|^2\|_{L^1(0,L_0)},  \quad k=1,2,\cdots
$$
if $\varsigma^\alpha |w_{\scriptscriptstyle0,\varsigma}|^2\in L^1(0,L_0)$ additionally, and
$$
\|w^{(k)}_0\|_{L^\infty(0,L_0)}\leq\|w_0\|_{L^\infty(0,L_0)}, \quad \|h_k \|_ {L^\infty(Q)} \leq \|h\|_ {L^\infty(Q)}, \quad k=1,2,\cdots
$$
if $w_0\in L^\infty(0,L_0)$ and $h\in L^\infty(Q)$ additionally.
Following the approach in \cite{Wang}, we consider the problem
\begin{equation}\label{2.3+}
\left\{\begin{array}{ll}
w^{(k)}_{t}-\big(a_k(\varsigma,t) w^{(k)}_\varsigma\big)_\varsigma- q_k(t)\varsigma w^{(k)}_\varsigma=h_k, & (\varsigma,t)\in Q,\\[2mm]
w^{(k)}(\text{\small$0$},t)=w^{(k)}(\text{\small $L_0$},t)=0, &t\in(0,T),\\[2mm]
w^{(k)}(\varsigma,\text{\small$0$})=w^{(k)}_0(\varsigma), &\varsigma\in (0,L_0).
\end{array}\right.
\end{equation}
According to the classical theory on parabolic equations, problem (\ref{2.3+})
admits a unique classical solution $w^{(k)}$. First, multiply the first equation of (\ref{2.3+}) by $w^{(k)}$ and integrate by parts over $Q_s$ $(0<s<T)$ to get
\begin{equation*}
\iint_{Q_s}\frac{1}{2}\big(|w^{(k)}|^2\big)_t+\big| a_k(\varsigma,t)w^{(k)}_\varsigma \big|^2+\frac{1}{2}q_k(t)|w^{(k)}|^2=\iint_{Q_s}h_kw^{(k)}.
\end{equation*}
Using Cauchy's inequality, we derive
\begin{equation*}
\iint_{Q_s}\frac{1}{2}\big(|w^{(k)}|^2\big)_t+|a_k(\varsigma,t)w^{(k)}_\varsigma|^2\leq\frac{1}{2}(\|q_k\|_{\infty}+1)\iint_{Q_s}|w^{(k)}|^2+\frac{1}{2}\iint_{Q_s}|h_k|^2,\quad \ 0<s<T.
\end{equation*}
Therefore, for any  $0<s<T$, it holds that
\begin{equation*}
\frac{1}{2}\int^{L_0}_{0}|w^{(k)}(s)|^2+\iint_{Q_s}|a_k(\varsigma,t)w^{(k)}_\varsigma|^2\leq\frac{1}{2}(\|q\|_{\infty}+1)\iint_{Q_s}|w^{(k)}|^2+\frac{1}{2}\big(\|h\|^2_{L^2(Q)}+\|w_0\|^2_{L^2(0,L_0)}\big).
\end{equation*}
Applying Gronwall's inequality, we obtain
\begin{equation*}
\|w^{(k)}\|^2_{L^\infty(0,T;L^2(0,L_0))}\leq e^{T(\|q\|_{\infty}+1)}\big(\|h\|^2_{L^2(Q)}+\|w_0\|^2_{L^2(0,L_0)}\big).
\end{equation*}
Furthermore,
\begin{equation}\label{e2.5-}
\|w^{(k)}\|^2_{L^\infty(0,T;L^2(0,L_0))}+\big\|a_k|w^{(k)}_\varsigma|^2\big\|_{L^1(Q)}\leq 2e^{2T(\|q\|_{\infty}+1)}\big(\|h_k\|^2_{L^2(Q)}+\|w_0^{(k)}\|^2_{L^2(0,L_0)}\big).
\end{equation}
From \eqref{e2.5-}, we deduce that there exists a function $w\in C([0,T];L^2(0,L_0))\cap L^2(0,T;H^1_\alpha(0,L_0))$ which is a weak solution to the problem \eqref{2.3} and satisfies the estimate $(i)$. Further estimates are established below.

Second, if $\varsigma^\alpha |w_{\scriptscriptstyle0,\varsigma}|^2\in L^1(0,L_0)$ additionally, multiplying the first equation of (\ref{2.3+}) by $w^{(k)}_t$ and integrating by parts over $Q_s$ $(0<s<T)$, we get
\begin{equation*}
\iint_{Q_s}|w^{(k)}_t|^2+\frac{1}{2}a_k(\varsigma,t)\big(|w^{(k)}_\varsigma|^2\big)_t-q_k(t)\varsigma w^{(k)}_\varsigma w^{(k)}_t=\iint_{Q_s}h_kw^{(k)}_t.
\end{equation*}
The above equation can be rewritten as
\begin{equation}\label{e2.5}
\iint_{Q_s}|w^{(k)}_t|^2+\frac{1}{2}\big(a_k(\varsigma,t)|w^{(k)}_\varsigma|^2\big)_t=\iint_{Q_s}q_k(t)\varsigma w^{(k)}_\varsigma w^{(k)}_t+h_kw^{(k)}_t+\frac{1}{2}\big(a_k(\varsigma,t)\big)_t|w^{(k)}_\varsigma|^2.
\end{equation}
Now, we estimate
\begin{equation}\label{2.9}
\begin{array}{ll}
&\quad\displaystyle\iint_{Q_s}q_k(t)\varsigma w^{(k)}_\varsigma w^{(k)}_t=\iint_{Q_s}\big(q_k(t){\textstyle L}^{\scriptstyle1-\frac{\alpha}{2}}_{\scriptscriptstyle0}\varsigma^{\frac{\alpha}{2}} w^{(k)}_\varsigma\big) \big({\textstyle L}^{\scriptstyle\frac{\alpha}{2}-1}_{\scriptscriptstyle0}\varsigma^{1-\frac{\alpha}{2}} w^{(k)}_t\big)\\[5mm]
&\leq\displaystyle\frac{1}{2}{{\textstyle L}^{\scriptstyle2-\alpha}_{\scriptscriptstyle0}}\|q\|^2_{\infty}\iint_{Q_s}\varsigma^\alpha |w^{(k)}_\varsigma|^2+\frac{1}{2}\iint_{Q_s}{{\textstyle L}^{\scriptstyle\alpha-2}_{\scriptscriptstyle0}}\varsigma^{2-\alpha}|w^{(k)}_t|^2\\[5mm]
&\leq\displaystyle\frac{1}{2}{{\textstyle L}^{\scriptstyle2-\alpha}_{\scriptscriptstyle0}}\|q\|^2_{\infty}\big\|{\textstyle\frac{1}{p}}\big\|_{\infty}\iint_{Q_s}a_k(\varsigma,t) |w^{(k)}_\varsigma|^2+\frac{1}{2}\iint_{Q_s}|w^{(k)}_t|^2.
\end{array}
\end{equation}
\eqref{e2.5}, together with \eqref{2.9}, indicates that
\begin{equation*}
\frac{1}{2}\iint_{Q_s}|w^{(k)}_t|^2+\big(a_k(\varsigma,t)|w^{(k)}_\varsigma|^2\big)_t\leq\iint_{Q_s}\frac{1}{2}\big({\textstyle L}^{\scriptscriptstyle2-\alpha}_{\scriptscriptstyle0}\|q\|^2_{\infty}\big\|{\textstyle\frac{1}{p}}\big\|_{\infty}+{\textstyle\frac{1}{a_k}\frac{\partial a_k}{\partial t}}\big)a_k(\varsigma,t)|w^{(k)}_\varsigma|^2+h_kw^{(k)}_t.
\end{equation*}
Since $\|\frac{1}{a_k}\frac{\partial a_k}{\partial t}\|_{L^\infty(Q)}\leq\|\frac{p'}{p}\|_{L^\infty(0,T)}$, it follows that
\begin{equation*}
\iint_{Q_s}\frac{1}{4}|w^{(k)}_t|^2+\frac{1}{2}\big(a_k(\varsigma,t)|w^{(k)}_\varsigma|^2\big)_t\leq\iint_{Q_s}\frac{1}{2}\big({\textstyle L}^{\scriptscriptstyle2-\alpha}_{\scriptscriptstyle0}\|q\|^2_{\infty}\big\|{\textstyle\frac{1}{p}}\big\|_{\infty}+{\textstyle\|\frac{p'}{p}\|_{\infty}}\big)a_k(\varsigma,t)|w^{(k)}_\varsigma|^2+\iint_{Q_s}h^2_k.
\end{equation*}
Therefore,
\begin{eqnarray*}
\frac{1}{4}\iint_{Q_s}|w^{(k)}_t|^2+\frac{1}{2}\int^{L_0}_0a_k(\varsigma,s)|w^{(k)}_\varsigma(s)|^2\leq\!\!\!\!\!\!\!\!&&\iint_{Q_s}\frac{1}{2}\big({\textstyle L}^{\scriptscriptstyle2-\alpha}_{\scriptscriptstyle0}\|q\|^2_{\infty}\big\|{\textstyle\frac{1}{p}}\big\|_{\infty}+{\textstyle\|\frac{p'}{p}\|_{\infty}}\big)a_k(\varsigma,t)|w^{(k)}_\varsigma|^2\\[2mm]
&&+\|h\|^2_{L^2(Q)}+\frac{1}{2}\|\varsigma^\alpha |w_{\scriptscriptstyle0,\varsigma}|^2\|_{L^1(0,L_0)}.
\end{eqnarray*}
Applying Gronwall's inequality, we obtain
\begin{eqnarray*}
\big\|a_k|w^{(k)}_\varsigma|^2\big\|_{L^\infty(0,T;L^1(Q))}
\leq e^{T({{\scriptstyle L}^{\scriptscriptstyle2-\alpha}_{\scriptscriptstyle0}}\|q\|^2_{\infty}\|\frac{1}{p}\|_{\infty}+\|\frac{p'}{p}\|_{\infty})}\big(2\|h\|^2_{L^2(Q)}+\|\varsigma^\alpha |w_{\scriptscriptstyle0,\varsigma}|^2\|_{L^1(0,L_0)}\big).
\end{eqnarray*}
Moreover,
\begin{equation}\label{e2.7}
\begin{array}{ll}
&\displaystyle\|w^{(k)}_t\|^2_{L^2(Q)}+\big\|a_k|w^{(k)}_\varsigma|^2\big\|_{L^\infty(0,T;L^1(Q))}\\[3mm]
\leq \!\!\!\!\!\!\!\!\!\!&\displaystyle4e^{2T({{\scriptstyle L}^{\scriptscriptstyle2-\alpha}_{\scriptscriptstyle0}}\|q\|^2_{\infty}\|\frac{1}{p}\|_{\infty}+\|\frac{p'}{p}\|_{\infty})}\big(2\|h\|^2_{L^2(Q)}+\|\varsigma^\alpha |w_{\scriptscriptstyle0,\varsigma}|^2\|_{L^1(0,L_0)}\big),
\end{array}
\end{equation}
which leads to the first estimate in $(ii)$.

Third, multiply (\ref{2.3+}) by $-\big(a_k(\varsigma,t) w^{(k)}_\varsigma\big)_\varsigma$ and then integrate
over $Q$ to get
\begin{equation*}
\iint_{Q}-w^{(k)}_t\big(a_k(\varsigma,t) w^{(k)}_\varsigma\big)_\varsigma+\big|\big(a_k(\varsigma,t) w^{(k)}_\varsigma\big)_\varsigma\big|^2+q_k(t)\varsigma w^{(k)}_\varsigma\big(a_k(\varsigma,t) w^{(k)}_\varsigma\big)_\varsigma=\iint_{Q}h_k\big(a_k(\varsigma,t) w^{(k)}_\varsigma\big)_\varsigma.
\end{equation*}
Rearranging the terms in the equation above, we arrive at
\begin{equation}\label{2.11}
\iint_{Q}\big|\big(a_k(\varsigma,t) w^{(k)}_\varsigma\big)_\varsigma\big|^2=\iint_{Q}\big(w^{(k)}_t+h_k\big)\big(a_k(\varsigma,t) w^{(k)}_\varsigma\big)_\varsigma-q_k(t)\varsigma w^{(k)}_\varsigma\big(a_k(\varsigma,t) w^{(k)}_\varsigma\big)_\varsigma.
\end{equation}
Notice that
$
\varsigma^{\alpha}=\frac{1}{p(t)}a(\varsigma,t)\leq\frac{1}{p(t)}a_k(\varsigma,t).
$
The last term of (\ref{2.11}) can be estimated as follows
\begin{equation}\label{2.12}
\begin{array}{ll}
&\displaystyle\iint_{Q}q_k(t)\varsigma w^{(k)}_\varsigma\big(a_k(\varsigma,t) w^{(k)}_\varsigma\big)_\varsigma=\iint_{Q}q_k(t)\varsigma^{1-\frac{\alpha}{2}}\big(\varsigma^{\frac{\alpha}{2}} w^{(k)}_\varsigma\big)\big(a_k(\varsigma,t) w^{(k)}_\varsigma\big)_\varsigma\\[5mm]
\leq\!\!\!\!\!\!\!\!\!\!&\displaystyle\iint_{Q}|q_k(t)|^2\varsigma^{2-\alpha}\big(\varsigma^{\alpha}|w^{(k)}_\varsigma|^2\big)+\frac{1}{4}\iint_{Q}\big|\big(a_k(\varsigma,t) w^{(k)}_\varsigma\big)_\varsigma\big|^2\\[5mm]
\leq\!\!\!\!\!\!\!\!\!\!&\displaystyle{\textstyle L}_{\scriptscriptstyle0}^{\scriptscriptstyle2-\alpha}\|q\|^2_{\infty}\big\|{\textstyle\frac{1}{p}}\big\|_{\infty}\iint_{Q}a_k(\varsigma,t)|w^{(k)}_\varsigma|^2+\frac{1}{4}\iint_{Q}\big|\big(a_k(\varsigma,t) w^{(k)}_\varsigma\big)_\varsigma\big|^2.
\end{array}
\end{equation}
Thus, we obtain by (\ref{2.11}) and (\ref{2.12})
\begin{equation*}
\iint_{Q}\frac{1}{4}\big|\big(a_k(\varsigma,t) w^{(k)}_\varsigma\big)_\varsigma\big|^2\leq{\textstyle L}_{\scriptscriptstyle0}^{\scriptscriptstyle2-\alpha}\|q\|^2_{\infty}\big\|{\textstyle\frac{1}{p}}\big\|_{\infty}\iint_{Q}a_k(\varsigma,t)|w^{(k)}_\varsigma|^2+\iint_{Q}|h_k|^2+|w^{(k)}_t|^2.
\end{equation*}
This, together with (\ref{e2.7}), indicates
\begin{eqnarray*}
&&\big\|\big(a_k w^{(k)}_\varsigma\big)_\varsigma\big\|^2_{L^2(Q)}\leq 16\big({\textstyle L}_{\scriptscriptstyle0}^{\scriptscriptstyle2-\alpha}\|q\|^2_{\infty}\big\|{\textstyle\frac{1}{p}}\big\|_{\infty}+1\big)\notag\\[2mm]
&&\quad\times e^{2T({{\scriptstyle L}^{\scriptscriptstyle2-\alpha}_{\scriptscriptstyle0}}\|q\|^2_{\infty}\|\frac{1}{p}\|_{\infty}+\|\frac{p'}{p}\|_{\infty})}
\big(3\|h\|^2_{L^2(Q)}+\|\varsigma^\alpha |w_{\scriptscriptstyle0,\varsigma}|^2\|_{L^1(0,L_0)}\big),
\end{eqnarray*}
which leads to the second estimate in $(ii)$.

Fourth, if $w_0\in L^\infty(0,L_0)$ and $h\in L^\infty(Q)$ additionally, then it follows from the maximum principle that
\begin{equation*}
\|w^{(k)}\|_{L^\infty(Q)}\leq T\|h_k\|_{L^\infty(Q)}+\|w_0^{(k)}\|_{L^\infty(Q)},
\end{equation*}
which yields the estimate in $(iii)$.

Finally, let us prove the uniqueness. Let $w_1$ and $w_2$ be two
weak solutions of the problem (\ref{2.3}) and denote
$$w(\varsigma,t)=w_1(\varsigma,t)-w_2(\varsigma,t),\quad (\varsigma,t)\in Q.$$
Then $w\in C([0,T];L^2(0,L_0))\cap L^2(0,T;H^1_\alpha(0,L_0))$ and for any function $\psi\in L^\infty(0,T;L^2(0,L_0))\cap L^2(0,T;H_\alpha^1(0,L_0))$ with $\psi_t\in L^2(Q)$ and $\psi(\cdot,T)|_{(0,L_0)}=0$, the following integral equality holds:
\begin{equation}\label{e2.9}
\iint_{Q}-w\psi_t+p(t)\varsigma^\alpha w_\varsigma \psi_\varsigma-q(t)\varsigma w_\varsigma \psi d\varsigma dt=0.
\end{equation}
For any $\zeta\in C^\infty_0(Q)$, the existence result shows that the problem
\begin{equation}\label{2.10}
\left\{\begin{array}{ll}
-\psi_{t}-p(t)(\varsigma^\alpha \psi_\varsigma)_\varsigma+ q(t)(\varsigma \psi)_\varsigma=\zeta, & (\varsigma,t)\in Q,\\[2mm]
\left\{\begin{array}{ll}
\psi(\text{\small$0$},t)=0, & \text{if}\quad \alpha\in[0,1),\\[2mm]
\lim\limits_{\varsigma\rightarrow0^+}(\varsigma^\alpha \psi_\varsigma)(\varsigma,t)=0, &\text{if}\quad \alpha\in[1,2),\\[2mm]
\end{array}\right.&\\[2mm]
\psi(\text{\small $L_0$},t)=0, &t\in(0,T),\\[2mm]
\psi(\varsigma,\text{\small$T$})=0, &\varsigma\in (0,L_0),
\end{array}\right.
\end{equation}
admits a weak solution $\psi\in C([0,T];L^2(0,L_0))\cap L^2(0,T;H^1_\alpha(0,L_0))$ with $\psi_t\in L^2(Q)$. Taking $\psi$ to be the solution of (\ref{2.10}) in \eqref{e2.9}, we get
\begin{equation*}
\iint_{Q}w\zeta=0.
\end{equation*}
This leads to
\begin{equation*}
w=0,\quad\ a.e.\ (\varsigma,t)\in Q,
\end{equation*}
owing to the arbitrariness of $\zeta\in C^\infty_0(Q)$. Therefore,
\begin{equation*}
w_1=w_2,\quad\ a.e.\ (\varsigma,t)\in Q.
\end{equation*}
The proof is complete.\endpf

By applying the inverse transformation $\Phi^{-1}_L$, the well-posedness result of (\ref{2.1}) (Theorem \ref{t2.1}) follows from Proposition \ref{p2.1}.

\subsection{Controllability problems and results}

Let us consider the controlled linear system
\begin{equation}\label{e2.1}
\left\{\begin{array}{ll}
y_t-(x^\alpha y_x)_x=\chi_\omega v, & (x,t)\in Q_L,\\[2mm]
\left\{\begin{array}{ll}
y(0,t)=0, & \text{if}\quad \alpha\in[0,1),\\[2mm]
\lim\limits_{x\rightarrow0}x^\alpha y_x(x,t)=0, &\text{if}\quad \alpha\in[1,2),\\[2mm]
\end{array}\right.&\\[2mm]
y(L(t),t)=0, &t\in(0,T),\\[2mm]
y(x,0)=y_0(x), &x\in (0,L_0).
\end{array}\right.
\end{equation}

Then we have the following approximate controllability result for (\ref{e2.1}).
\begin{theorem}\label{c2.1}
For any $y_0\in L^2(0, L_0)$ and any $\beta>0$, there exist pairs $(v_\beta, y_\beta )$, with
\begin{equation*}
v_\beta\in L^2(\omega\times(0,T)),\quad  y_\beta\in C([0,T];L^2(0,L(t)))\cap L^2(0,T;H_\alpha^1(0,L(t))),
\end{equation*}
satisfying (\ref{e2.1}) and
\begin{eqnarray*}\label{bet1}
\|y_\beta(\cdot,T )\|_{L^2(0,L(T))}\leq\beta.
\end{eqnarray*}
Furthermore, $v_\beta$ can be found such that
\begin{eqnarray}\label{cost1}
\|v_\beta\|_{L^2(\omega\times(0,T))}\leq C\|y_0\|_{L^2(0,L_0)},
\end{eqnarray}
where $C$ is a positive constant depending on $\alpha$, $N_L, L_\ast, B,$ $\omega$ and $T$, but independent of $\beta$.
\end{theorem}

An immediate consequence of Theorem \ref{c2.1} is the following null controllability of (\ref{e2.1}):
\begin{corollary}\label{cor2.1}
For any $y_0\in L^2(0, L_0)$, there exist pairs $(v, y)$, with
\begin{equation*}
v\in L^2(\omega\times(0,T)),\quad y\in  C([0,T];L^2(0,L(t)))\cap L^2(0,T;H_\alpha^1(0,L(t))),
\end{equation*}
satisfying (\ref{e2.1}) and
\begin{eqnarray*}\label{bet2}
y(x,T )=0,\quad\quad x\in(0,L(T)).
\end{eqnarray*}
Furthermore, $v$ can be found such that
\begin{eqnarray}\label{cost2}
\|v\|_{L^2(\omega\times(0,T))}\leq C\|y_0\|_{L^2(0,L_0)},
\end{eqnarray}
where $C$ is a positive constant depending only on $\alpha$, $N_L, L_\ast, B,$ $\omega$ and $T$.
\end{corollary}

\subsection{A Carleman estimate}
The proof of Theorem \ref{c2.1} relies on duality arguments. The main tool is a global Carleman
estimate for the adjoint system of (\ref{e2.1}):
\begin{equation}\label{2.8}
\left\{\begin{array}{ll}
\varphi_{t}+(x^\alpha \varphi_x)_x=g, & (x,t)\in Q_L,\\[2mm]
\left\{\begin{array}{ll}
\varphi(0,t)=0, & \text{if}\quad \alpha\in[0,1),\\[2mm]
\lim\limits_{x\rightarrow0^+}x^\alpha \varphi_x(x,t)=0, &\text{if}\quad \alpha\in[1,2),\\[2mm]
\end{array}\right. &\\[2mm]
\varphi(L(t),t)=0, &t\in(0,T),\\[2mm]
\varphi(x,T)=\varphi_T(x), &x\in (0, L(T)),
\end{array}\right.
\end{equation}
where $g \in L^2(Q_L)$ and $\varphi_T\in L^2(0,L(T)).$

As a preliminary, we introduce the following Hardy-type inequalities in non-cylindrical domains as a major ingredient in studying degenerate problems (for the reader's
convenience, we present the proof in Section \ref{sec6}).
\begin{lemma}\label{lht}(Hardy-type inequalities)
$(i)$ Let $0 \leq \alpha^\star< 1$. If $z$ is locally absolutely continuous on
$(0, L(t))$ for almost
every $t\in(0, T)$, and satisfies
\begin{eqnarray*}
z(x) \rightarrow0,\ x\rightarrow0^+ \ \text{and} \ \int^{L(t)}_0 x^{\alpha^\star} z^2_x < \infty,
\end{eqnarray*}
then the following inequality holds:
\begin{eqnarray}\label{ht}
\int^{L(t)}_0x^{\alpha^\star-2}z^2\leq \frac{4}{(\alpha^\star-1)^2}\int^{L(t)}_0 x^{\alpha^\star}z_x^2.
\end{eqnarray}
$(ii)$ Let $1 <\alpha^\star< 2$. Then the above inequality (\ref{ht}) still holds if $z$ is locally absolutely continuous on
$(0, L(t))$ for almost
every $t\in(0, T)$, and satisfies
\begin{eqnarray*}
z(x)\rightarrow0,\ x\rightarrow L(t)^- \ \text{and} \ \int^{L(t)}_0 x^{\alpha^\star} z^2_x < \infty.
\end{eqnarray*}
\end{lemma}

We have the following global Carleman estimate for (\ref{2.8}).
\begin{theorem}\label{t2.2}
Let $0\leq\alpha<2$ and $T>0$ be given. Then there exist a function $\sigma: [0,B]\times(0, T)\rightarrow \mathbb R_+$ of the form $\sigma(x,t) = \theta(t)\Phi(x)$, with
\begin{equation}\label{2.20}
\Phi(x)>0\ \  \forall x\in[0,B] \quad \ \text{and}\quad \ \theta(t)\rightarrow \infty \ \ \text{as}\ \ t \rightarrow 0^+,\ T^-,
\end{equation}
and two positive constants $s_0$ and $C$, depending only on $\alpha$, $N_L$, $L_\ast$,
$B$, $\omega$ and $T$, such that for any $s\geq s_0$, any $g\in
L^2(Q_L)$ and any $\varphi_T\in L^2(0,L(T))$, the solution $\varphi$ of (\ref{2.8}) satisfies
\begin{eqnarray*}
\iint_{Q_L}\big(s\theta x^\alpha \varphi^2_x+s^3\theta^3x^{2-\alpha}\varphi^2\big)e^{-2s\sigma}dxdt+s\int^T_0L^\alpha(t) e^{-2s\sigma(L(t),t)}\varphi^2_x(L(t),t)dt\\[2mm]
\leq C\(\int^T_0\int_{\omega_0}\big(s\theta x^\alpha \varphi^2_x+s^3\theta^3x^{2-\alpha}\varphi^2\big)e^{-2s\sigma}dxdt+\iint_{Q_L}g^2e^{-2s\sigma}dxdt\).
\end{eqnarray*}
\end{theorem}

The proof will be presented in Section \ref{carle}.

\subsection{Observability inequality}
As an application of Theorem \ref{t2.2}, we will derive an observability inequality for the homogeneous parabolic system:
\begin{equation}\label{bf}
\left\{\begin{array}{ll}
\varphi_t+(x^\alpha \varphi_x)_x=0, & (x,t)\in Q_L,\\[2mm]
\left\{\begin{array}{ll}
\varphi(0,t)=0, & \text{if}\quad \alpha\in[0,1),\\[2mm]
\lim\limits_{x\rightarrow0^+}x^\alpha \varphi_x(x,t)=0, &\text{if}\quad \alpha\in[1,2),\\[2mm]
\end{array}\right.&\\[2mm]
\varphi(L(t),t)=0, &t\in(0,T),\\[2mm]
\varphi(x,T)=\varphi_T(x), &x\in (0,L(T)),
\end{array}\right.
\end{equation}
where $\varphi_T\in L^2(0, L(T ))$. Our result is the following.
\begin{proposition}\label{pr2.1}
There exists a constant $C > 0$, depending only on $\alpha$, $N_L$, $L_\ast$,
$B$, $\omega$ and $T$, such that for any $\varphi_T\in L^2(0, L(T ))$, the associated
solution to (\ref{bf}) satisfies
\begin{equation}\label{ob}
 \int^{L_0}
_0
|\varphi(x, 0)|^2
dx \leq C \iint_{\omega\times(0,T )}|\varphi(x, t)|^2dxdt.
\end{equation}
\end{proposition}

Before proving Proposition \ref{pr2.1}, we provide the Caccioppoli-type inequality for (\ref{bf}), whose proof is given in Section \ref{sec6}.
\begin{lemma}\label{cacc}(Caccioppoli-type inequality) For all $s>0$, the solution $\varphi$ of (\ref{bf}) satisfies
\begin{equation*}
\int^T_0\int_{\omega_0} e^{-2s\sigma}\varphi_x^2dxdt\leq C(s,T)\int^T_0\int_{\omega}\varphi^2dxdt.
\end{equation*}
\end{lemma}

Combing Theorem \ref{t2.2} with Lemma \ref{cacc}, one can easily get the following Carleman estimate for (\ref{bf}).
\begin{lemma}\label{l2.2}
Let $0\leq\alpha<2$ and $T>0$ be given. Then there exist a function $\sigma$ of the form \eqref{2.20}
and a positive constant $s_0$, depending only on $\alpha$, $N_L$, $L_\ast$,
$B$, $\omega$ and $T$, such that for any $s\geq s_0$, one can find a positive constant $C$ so that for any $\varphi_T\in L^2(0,L(T))$, the solution $\varphi$ of (\ref{bf}) satisfies
\begin{eqnarray*}
\iint_{Q_L}\big(s\theta x^\alpha \varphi^2_x+s^3\theta^3x^{2-\alpha}\varphi^2\big)e^{-2s\sigma}dxdt+s\int^T_0L^\alpha(t) e^{-2s\sigma(L(t),t)}\varphi^2_x(L(t),t)dt
\leq C\int^T_0\int_{\omega}\varphi^2dxdt.
\end{eqnarray*}
\end{lemma}

Now, we are in a position to present the proof of Proposition \ref{pr2.1}. Note that the moving boundary alters the weighted energy integral $\int_0^{L(t)}x^\alpha\varphi^2_x(x,t)dx$ for degenerate parabolic equations, which is different from the case of a fixed cylindrical domain studied in \cite{4}. Therefore, to establish the observability inequality, we need to incorporate the $L^2$-norm estimate for solutions to (\ref{bf}).

\emph{Proof of Proposition \ref{pr2.1}.} The proof proceeds in two steps.

{\bf Step 1.} First, multiply both sides of the first equation in (\ref{bf}) by $\varphi_t$ and integrate over $(0, L(t))$ to get
\begin{equation}\label{int}
\begin{array}{ll}
0\!\!\!&\displaystyle=\int_0^{L(t)}\big(\varphi_t+(x^\alpha\varphi_x)_x\big)\varphi_tdx\\[4.5mm]
&\displaystyle=\int_0^{L(t)}\varphi^2_tdx+\Big[x^\alpha\varphi_x\varphi_{t}\Big]^{L(t)}_0-\int_0^{L(t)}x^\alpha\varphi_x\varphi_{tx}dx\quad\ \forall t\in(0,T).
\end{array}
\end{equation}
Next, we analyze the last two terms in \eqref{int}. Using the boundary condition $\varphi(L(t),t)=0$ (which implies $\varphi_x(L(t),t)L'(t)+\varphi_t(L(t),t)=0$), we derive
\begin{eqnarray}\label{int1}
\big[x^\alpha\varphi_x\varphi_{t}\big]^{L(t)}_0=-L^\alpha(t)L'(t)\varphi^2_x(L(t),t).
\end{eqnarray}
On the other hand, it is easy to see that
\begin{eqnarray}\label{int2}
\int_0^{L(t)}x^\alpha\varphi_x\varphi_{tx}dx=\frac{d}{dt}\bigg(\int_0^{L(t)}\frac{1}{2}x^\alpha\varphi^2_xdx\bigg)-\frac{1}{2}L'(t)L^\alpha(t)\varphi^2_x(L(t),t).
\end{eqnarray}
Substituting (\ref{int1}) and (\ref{int2}) into (\ref{int}) yields
\begin{eqnarray*}
0\geq-\frac{d}{dt}\bigg(\int_0^{L(t)}\frac{1}{2}x^\alpha\varphi^2_xdx\bigg)-\frac{1}{2}L^\alpha(t)L'(t)\varphi^2_x(L(t),t).
\end{eqnarray*}
This implies that for any $t\in(T/4,T)$,
\begin{eqnarray*}
\int_0^{L(T/4)}x^\alpha\varphi^2_x(x,T/4)dx\leq\int_0^{L(t)}x^\alpha\varphi^2_x(x,t)dx+\int^t_{T/4}L^\alpha(\tau)L'(\tau)\varphi^2_x(L(\tau),\tau)d\tau.
\end{eqnarray*}
Integrating the above inequality over $(T/4, 3T/4)$, we obtain
\begin{eqnarray*}
&&\int_0^{L(T/4)}x^\alpha\varphi^2_x(x,T/4)dx\leq\frac{2}{T}\int^{3T/4}_{T/4}\bigg(\int_0^{L(t)}x^\alpha\varphi^2_x(x,t)dx+\int^t_{T/4}L^\alpha(\tau)L'(\tau)\varphi^2_x(L(\tau),\tau)d\tau\bigg)dt\\[2.5mm]
&&\quad\quad\leq C(s,T)\bigg(\int^{3T/4}_{T/4}\int_0^{L(t)}\theta x^\alpha\varphi^2_x(x,t)e^{-2s\sigma}dxdt+N_L\int^{3T/4}_{T/4}L^\alpha(t)e^{-2s\sigma(L(t),t)}\varphi^2_x(L(t),t)dt\bigg).
\end{eqnarray*}
Hence, owing to Lemma \ref{l2.2},
\begin{eqnarray}\label{e1.13}
\int_0^{L(T/4)}x^\alpha\varphi^2_x(x,T/4)dx\leq C\int^T_0\int_{\omega}\varphi^2dxdt.
\end{eqnarray}
For $\alpha \neq1$, applying Hardy's inequality (Lemma \ref{lht}) with
$\alpha^\star = \alpha$, we deduce from \eqref{e1.13} that
\begin{equation}\label{ec1}
\int_0^{L(T/4)}x^{\alpha-2}\varphi^2(x,T/4)dx\leq C\int_0^{L(T/4)}x^\alpha\varphi^2_x(x,T/4)dx\leq C\int^T_0\int_{\omega}\varphi^2dxdt.
\end{equation}
In the case of $\alpha =1$, from \eqref{e1.13}, we deduce that, for any $0<\delta<1$,
\begin{equation*}
\int_0^{L(T/4)}x^{1+\delta}\varphi^2_x(x,T/4)dx\leq \int_0^{L(T/4)}xL(T/4)^{\delta}\varphi^2_x(x,T/4)dx\leq C\int^T_0\int_{\omega}\varphi^2dxdt.
\end{equation*}
Now, applying Hardy's inequality with $\alpha^\star=1+ \delta$, we obtain
\begin{equation}\label{ec2}
\int_0^{L(T/4)}x^{\delta-1}\varphi^2(x,T/4)dx\leq C\int_0^{L(T/4)}x^{1+\delta}\varphi^2_x(x,T/4)dx\leq C\int^T_0\int_{\omega}\varphi^2dxdt.
\end{equation}
In both cases, notice that $\alpha <2$ and $L_\ast\leq L(t)\leq B$ for all $t\in[0, T ]$. It follows  from \eqref{ec1} or \eqref{ec2} that
\begin{eqnarray}\label{ob1}
\int_0^{L(T/4)}\varphi^2(x,T/4)dx\leq C\int^T_0\int_{\omega}\varphi^2dxdt.
\end{eqnarray}

{\bf Step 2.} Multiplying both sides of the first equation in (\ref{bf}) by $\varphi$ and integrating over $(0,L(t))$, we get
\begin{eqnarray*}
0=\int_0^{L(t)}\big(\varphi_t+(x^\alpha\varphi_x)_x\big)\varphi dx=\frac{1}{2}\frac{d}{dt}\bigg(\int_0^{L(t)}\varphi^2dx\bigg)-\int_0^{L(t)}x^\alpha\varphi^2_xdx\quad\ \forall t\in(0,T).
\end{eqnarray*}
Hence,
\begin{eqnarray*}
\frac{d}{dt}\bigg(\int_0^{L(t)}\varphi^2dx\bigg)\geq0.
\end{eqnarray*}
This implies
\begin{eqnarray}\label{var2}
\int_0^{L_0}\varphi^2(x,0)dx\leq\int_0^{L(T/4)}\varphi^2(x,T/4)dx.
\end{eqnarray}
Combining (\ref{ob1}) with (\ref{var2}), we obtain the desired inequality (\ref{ob}).\endpf

\subsection{Controllability analysis for (\ref{e2.1})}\label{s2.5}

By duality and by applying the variational approach, the approximate controllability result in Theorem \ref{c2.1} can be established. To this
end, for any $y_0\in L^2(0, L_0)$ and given $\beta > 0$, let us define the functional $J_{\beta}(\cdot,L)$:
\begin{equation}\label{j}
J_{\beta}(\varphi_T,L):=\frac{1}{2}\iint_{\omega\times(0,T )}|\varphi|^2dxdt
+\beta\|\varphi_T\|_{ L^2(0, L(T ))} +\big(\varphi(\cdot, 0), y_0\big)_{ L^2(0, L_0)},
\end{equation}
where $\varphi$ is the solution of the adjoint system (\ref{bf}) with terminal data $\varphi_T$ at time $t = T$.

We shall show that $J_{\beta}(\cdot,L)$ is continuous, strictly convex, and coercive in $L^2(0, L(T ))$. These three properties ensure that $J_{\beta}(\cdot,L)$ possesses a unique minimizer $\hat{\varphi}_T$.

$(1)$ From Theorem \ref{t2.1} $(i)$, it is obvious that $J_{\beta}(\cdot,L)$ is continuous linear functional from $L^2(0, L(T ))$ to $\mathbb R$.

$(2)$ We show that $J_{\beta}(\cdot,L)$ is
strictly convex.

Let $\varphi^0_T,\varphi^1_T\in L^2(0, L(T ))$ and $\lambda\in(0, 1)$. We have
\begin{equation*}
J_{\beta}\big(\lambda\varphi^0_T+(1-\lambda)\varphi^1_T,L\big)\leq\lambda J_{\beta}\big(\varphi^0_T,L)+(1-\lambda)J_{\beta}\big(\varphi^1_T,L\big)-\frac{\lambda(1-\lambda)}{2}\iint_{\omega\times(0,T )}|\varphi^0-\varphi^1|^2dxdt,
\end{equation*}
where $\varphi^0$ and $\varphi^1$  are the solutions to (\ref{bf}) with terminal data $\varphi^0_T$ and $\varphi^1_T$, respectively.

From Proposition \ref{pr2.1}, there exists a positive constant $C$ (only depending on $\alpha, N_L, L_\ast,
B, \omega$ and $T$) such that
\begin{eqnarray*}
\iint_{\omega\times(0,T )}|\varphi^0-\varphi^1|^2dxdt\geq C\|\varphi^0(\cdot, 0)-\varphi^1(\cdot, 0)\|_{L^2(0, L_0)}.
\end{eqnarray*}

Consequently, for any $\varphi^0_T\neq\varphi^1_T$,
\begin{eqnarray*}
J_{\beta}\big(\lambda\varphi^0_T+(1-\lambda)\varphi^1_T,L\big)<\lambda J_{\beta}\big(\varphi^0_T,L)+(1-\lambda)J_{\beta}\big(\varphi^1_T,L\big)
\end{eqnarray*}
and $J_{\beta}(\cdot,L)$ is strictly convex.

$(3)$ We prove that the functional $J_{\beta}(\cdot,L)$ is also coercive, i.e.
\begin{eqnarray}\label{2.29}
\lim\limits_{\|\varphi_T\|_{L^2(0, L(T ))}\rightarrow\infty}J_{\beta}(\varphi_T,L)=\infty.
\end{eqnarray}
In fact, we shall prove that
\begin{eqnarray}\label{2.30}
\lim\limits_{\|\varphi_T\|_{L^2(0, L(T ))}\rightarrow\infty}\inf\frac{ J_{\beta}(\varphi_T,L)}{\|\varphi_T\|_{L^2(0, L(T ))}}\geq\varepsilon.
\end{eqnarray}
Evidently, (\ref{2.30}) implies (\ref{2.29}). In order to prove (\ref{2.30}), let $\{\varphi_{T,j}\}^\infty_{j=1}\subset L^2(0, L(T ))$ be a sequence of terminal data for
the adjoint system (\ref{bf})  with $\|\varphi_{T,j}\|_{L^2(0, L(T ))}\rightarrow\infty$. We normalize them
\begin{eqnarray*}
\tilde{\varphi}_{T,j}=\frac{\varphi_{T,j}}{\|\varphi_{T,j}\|},
\end{eqnarray*}
so that $\|\tilde{\varphi}_{T,j}\|_{L^2(0, L(T ))}=1$.

On the other hand, let $\tilde{\varphi}_{j}$ be the solution of (\ref{bf}) with terminal data $\tilde{\varphi}_{T,j}$. Then
\begin{eqnarray*}
\frac{J_{\beta}(\varphi_{T,j},L)}{\|\varphi_{T,j}\|_{L^2(0, L(T ))}}=\frac{1}{2}\|\varphi_{T,j}\|_{L^2(0, L(T ))}\iint_{\omega\times(0,T )}|\tilde{\varphi}_j|^2dxdt
+\beta +\big(\tilde{\varphi}_j(\cdot, 0), y_0\big)_{ L^2(0, L_0)}.
\end{eqnarray*}
The following two cases may occur:

a) $\lim\limits_{j\rightarrow\infty}\inf\iint_{\omega\times(0,T )}|\tilde{\varphi}_j|^2dxdt>0$. In this case we obtain immediately that
\begin{eqnarray*}
\frac{J_{\beta}(\varphi_{T,j},L)}{\|\varphi_{T,j}\|_{L^2(0, L(T ))}}\rightarrow\infty.
\end{eqnarray*}

b) $\lim\limits_{j\rightarrow\infty}\inf\iint_{\omega\times(0,T )}|\tilde{\varphi}_j|^2dxdt=0$.  In this case since $\|\tilde{\varphi}_{T,j}\|_{L^2(0, L(T ))}=1$, we can  extract a subsequence (still denoted by $\tilde{\varphi}_{T,j}$), which weakly converges to an element $\tilde{\varphi}_{T}$ in $L^2(0, L(T ))$. From Theorem \ref{t2.1} $(i)$, $\tilde{\varphi}_{j}$ converges weakly in $L^2(Q_L)$ to the solution  $\tilde{\varphi}$ of (\ref{bf}) with $\tilde{\varphi}(T)=\tilde{\varphi}_{T}$. Moreover, by weakly
lower semi-continuity,
\begin{eqnarray*}
\iint_{\omega\times(0,T )}|\tilde{\varphi}|^2dxdt\leq\lim\limits_{j\rightarrow\infty}\inf\iint_{\omega\times(0,T )}|\tilde{\varphi}_j|^2dxdt=0
\end{eqnarray*}
and therefore $\tilde{\varphi}=0$ in $\omega\times(0,T )$.

Lemma \ref{l2.2} implies that  $\tilde{\varphi}=0$ in $Q_L$ and
consequently $\tilde{\varphi}_{T}=0$.

Therefore, $\tilde{\varphi}_{T,j}\rightharpoonup0$ in $L^2(0, L(T ))$ and consequently $\tilde{\varphi}_j(\cdot, 0)\rightharpoonup0$ in $L^2(0, L_0)$.

Hence
\begin{eqnarray*}
\lim\limits_{j\rightarrow\infty}\inf\frac{J_{\beta}(\varphi_{T,j},L)}{\|\varphi_{T,j}\|_{L^2(0, L(T ))}}\geq \lim\limits_{j\rightarrow\infty}\inf \big[\beta+\big(\tilde{\varphi}_j(\cdot, 0), y_0\big)_{ L^2(0, L_0)}\big]=\beta,
\end{eqnarray*}
and (\ref{2.30}) follows.

Denote by $\hat{\varphi}$ the corresponding solution to (\ref{bf}) associated with $\hat{\varphi}_T$. Then, the restriction
\begin{equation}\label{con}
v_\beta = \hat{\varphi}|_{\omega\times(0,T )}
\end{equation}
constitutes an approximate control for system (\ref{e2.1}) with initial data $y_0$.

Indeed, if $J_{\beta}(\cdot,L)$ attains its minimum at $\hat{\varphi}_T\in L^2(0, L(T ))$, then for any $\psi_T\in L^2(0, L(T ))$ and $h\in \mathbb R$, we have $J_{\beta}(\hat{\varphi}_T,L)\leq J_{\beta}(\hat{\varphi}_T+h\psi_T,L)$. By the definition of $J_{\beta}(\cdot,L)$ in (\ref{j}),
\begin{eqnarray*}
&& J_{\beta}(\hat{\varphi}_T + h\psi_T, L) \\[2mm]
&& = \frac{1}{2} \iint_{\omega \times (0,T)} |\hat{\varphi} + h\psi|^2 dxdt
 + \beta \|\hat{\varphi}_T + h\psi_T\|_{L^2(0, L(T))}
 + \left( \hat{\varphi}(\cdot, 0) + h\psi(\cdot, 0), y_0 \right)_{L^2(0, L_0)}\\[2mm]
&& = \frac{1}{2} \iint_{\omega \times (0,T)} |\hat{\varphi}|^2 \, dxdt
 + \frac{h^2}{2} \iint_{\omega \times (0,T)} |\psi|^2  dxdt
 + h \iint_{\omega \times (0,T)} \hat{\varphi}\psi  dxdt \\[2mm]
&& \quad + \beta \|\hat{\varphi}_T + h\psi_T\|_{L^2(0, L(T))}
 + \left( \hat{\varphi}(\cdot, 0) + h\psi(\cdot, 0), y_0 \right)_{L^2(0, L_0)}.
\end{eqnarray*}
Thus
\begin{eqnarray*}
0\!\!\!\!\!\!\!\!&&\leq J_{\beta}(\hat{\varphi}_T+h\psi_T,L)-J_{\beta}(\hat{\varphi}_T,L)\\[2mm]
\!\!\!\!\!\!\!\!&&=\frac{h^2}{2}\iint_{\omega\times(0,T )}|\psi|^2dxdt+h\iint_{\omega\times(0,T )}\hat{\varphi}\psi dxdt\\[2mm]
&&\quad+\beta\|\hat{\varphi}_T+h\psi_T\|_{ L^2(0, L(T ))} -\beta\|\hat{\varphi}_T\|_{ L^2(0, L(T ))}+h\big(\psi(\cdot, 0), y_0\big)_{ L^2(0, L_0)},
\end{eqnarray*}
Since
\begin{eqnarray*}
\|\hat{\varphi}_T+h\psi_T\|_{ L^2(0, L(T ))}-\|\hat{\varphi}_T\|_{ L^2(0, L(T ))}\leq|h|\|\psi_T\|_{ L^2(0, L(T ))},
\end{eqnarray*}
we arrive at
\begin{eqnarray*}
0\!\!\!\!\!\!\!\!&&\leq J_{\beta}(\hat{\varphi}_T+h\psi_T,L)-J_{\beta}(\hat{\varphi}_T,L)\\[2mm]
\!\!\!\!\!\!\!\!&&\leq\frac{h^2}{2}\iint_{\omega\times(0,T )}|\psi|^2dxdt+h\iint_{\omega\times(0,T )}\hat{\varphi}\psi dxdt\\[2mm]
&&\quad+\beta|h|\|\psi_T\|_{ L^2(0, L(T ))}+h\big(\psi(\cdot, 0), y_0\big)_{ L^2(0, L_0)},
\end{eqnarray*}
for all $h\in \mathbb R$ and $\psi_T\in L^2(0, L(T ))$.

Dividing the above inequality by $h>0$ and taking the limit as $h\rightarrow0$, we obtain
\begin{eqnarray*}
0\leq\iint_{\omega\times(0,T )}\hat{\varphi}\psi dxdt+\beta\|\psi_T\|_{ L^2(0, L(T ))}+\big(\psi(\cdot, 0), y_0\big)_{ L^2(0, L_0)}.
\end{eqnarray*}
Performing the same calculations for $h < 0$ gives
\begin{eqnarray*}
0\geq\iint_{\omega\times(0,T )}\hat{\varphi}\psi dxdt-\beta\|\psi_T\|_{ L^2(0, L(T ))}+\big(\psi(\cdot, 0), y_0\big)_{ L^2(0, L_0)}.
\end{eqnarray*}
Combining these two results, we conclude that
\begin{eqnarray}\label{2.19}
\left|\iint_{\omega\times(0,T )}\hat{\varphi}\psi dxdt+\big(\psi(\cdot, 0), y_0\big)_{ L^2(0, L_0)}\right|\leq\beta\|\psi_T\|_{ L^2(0, L(T ))},\quad \forall \psi_T\in L^2(0, L(T )).
\end{eqnarray}
On the other hand, if we take the control $f=\hat{\varphi}$ in \eqref{e2.1}, then by multiplying the first equation of \eqref{e2.1} by $\psi$ (which is the solution to (\ref{bf}) with terminal data $\psi_T$) and integrating by parts over $Q_L$, we get
\begin{eqnarray*}
\iint_{\omega\times(0,T )}\hat{\varphi}\psi dxdt=\big(\psi_T, y(\cdot,T)\big)_{ L^2(0, L(T))}-\big(\psi(\cdot, 0), y_0\big)_{ L^2(0, L_0)}.
\end{eqnarray*}
Substituting the above formula into \eqref{2.19}, we find
\begin{eqnarray*}
\left|\big(\psi_T, y(\cdot,T)\big)_{ L^2(0, L(T))}\right|\leq\beta\|\psi_T\|_{ L^2(0, L(T ))},\quad \forall \psi_T\in L^2(0, L(T )),
\end{eqnarray*}
which is equivalent to
\begin{eqnarray*}
\|y(\cdot,T)\|_{ L^2(0, L(T ))}\leq\beta.
\end{eqnarray*}
This implies that the minimizer of $J_{\beta}(\cdot,L)$ provides an approximate control of the form \eqref{con} for problem \eqref{e2.1}.

Finally, we establish the uniform cost estimate for the control function $v_\beta$.
From the inequality $J_{\beta}(\hat{\varphi}_{\scriptscriptstyle T},L)\leq J_{\beta}(0,L) = 0$ and
(\ref{ob}), we deduce that
\begin{eqnarray*}
&&\frac{1}{2}\iint_{\omega\times(0,T )}|\hat{\varphi}|^2dxdt
+\beta\|\hat{\varphi}_{\scriptscriptstyle T}\|_{ L^2(0, L(T ))}
 \leq-\big(\hat{\varphi}(\cdot, 0), y_0\big)_{ L^2(0, L_0)}\\[2mm]
\leq\!\!\!\!\!\!\!\!&&\frac{1}{4C}\|\hat{\varphi}(\cdot, 0)\|^2_{ L^2(0, L_0)}+ C\|y_0\|^2_{L^2(0,L_0)}
\leq
\frac{1}{4} \iint_{\omega\times(0,T )}|\hat{\varphi}|^2dxdt+ C\|y_0\|^2_{L^2(0,L_0)}.
\end{eqnarray*}
Therefore,
\begin{eqnarray}\label{2.26}
&&\iint_{\omega\times(0,T )}|\hat{\varphi}|^2dxdt
+4\beta\|\hat{\varphi}_{\scriptscriptstyle T}\|_{ L^2(0, L(T ))}
\leq 4C\|y_0\|^2_{L^2(0,L_0)}.
\end{eqnarray}
Consequently,
\begin{eqnarray*}
\|v_\beta\|^2_{L^2(\omega\times(0,T ))} \leq C\|y_0\|^2_{L^2(0,L_0)},
\end{eqnarray*}
where $C$ is a positive constant depending only on $\alpha$, $N_L, L_\ast, B,$ $\omega$ and $T$.

\subsection{A regularity property}\label{sec3.1}

In this part, we establish the H\"{o}lder regularity for solutions to degenerate parabolic equations in non-cylindrical domains. This regularity result will play a crucial role in the proof of Theorem 1.1 in Section \ref{sec3}.

Let $(v, y)$ be the control-state pair furnished by Theorem \ref{c2.1} or Corollary \ref{cor2.1} and let us introduce the set $R_L:= Q_L \cap \{(x, t) : x > b'\},$ where $b < b'< L_\ast$. Following an argument analogous to that presented in \cite{FLde}, we conclude that there exists a positive constant
$\kappa \in (0, 1/2]$ depending only on $\alpha$, $N_L$, $L_\ast$, $B$, $\omega$ and $T$ such that $$y \in C_{x,t}^{1+\kappa,\kappa/2}(\overline{R_L}),$$ where $C_{x,t}^{1+\kappa,\kappa/2}(\overline{R_L})$ is the space of functions $z\in C(\overline{R_L})$ that possess a continuous partial derivative with respect to $x$ in $\overline{R_L}$ and satisfy
\begin{eqnarray*}
\sup\limits_
{(x,t),(x',t')\in\overline{R_L}}\bigg(\frac{|z(x, t) -z(x',t')|}{
|x - x'| +|t -t'|^{\kappa/2}}+\frac{|z_x(x,t)-z_x(x',t')|}{
|x- x'|^{\kappa}+|t-t'|^{\kappa/2}}\bigg)<+\infty.
\end{eqnarray*}
\begin{figure}[htbp]
  \centering
  \includegraphics[width=0.36\textwidth]{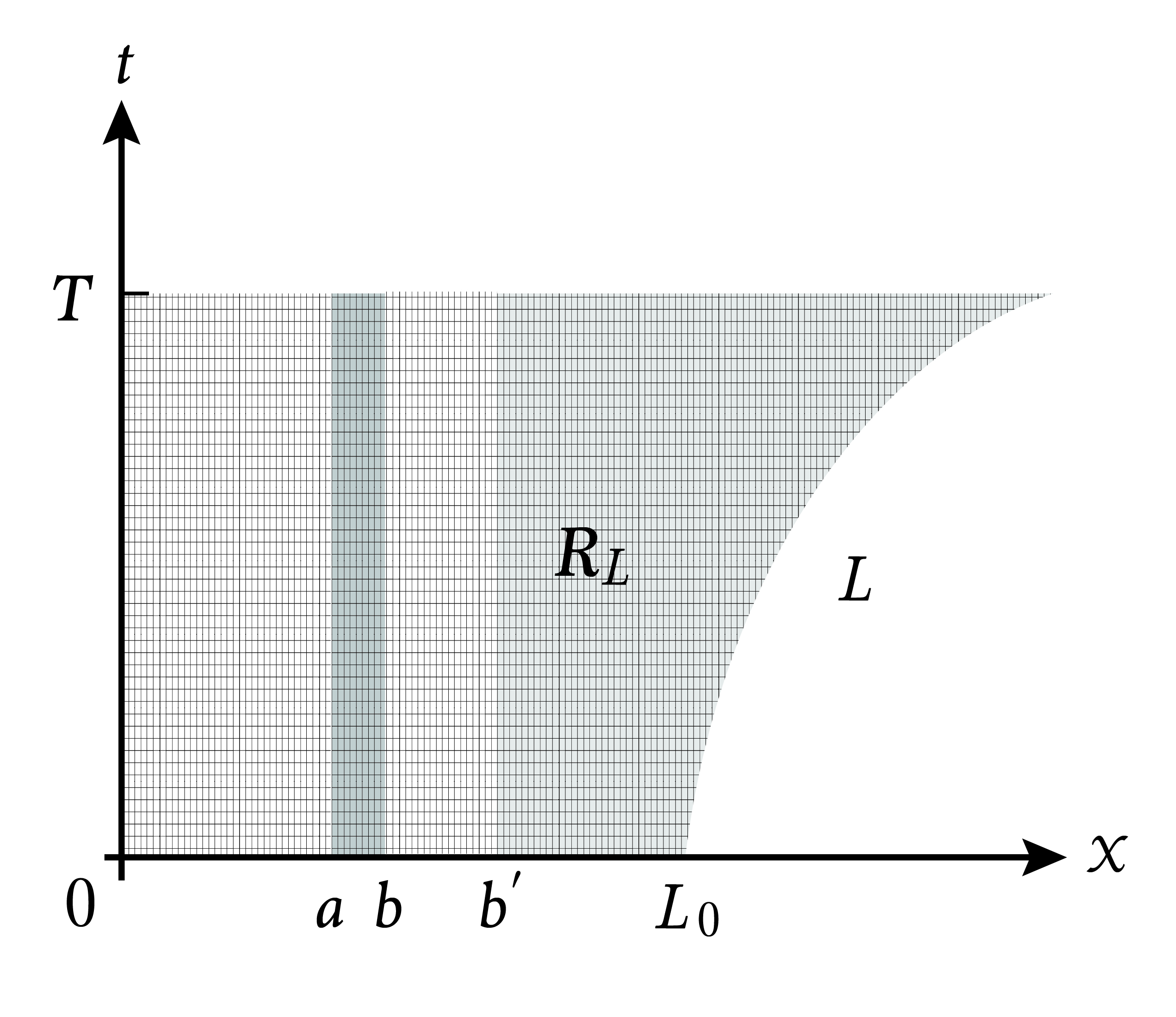}\\
  \caption{Diagram of region $R_L$}\label{Image}
\end{figure}
Now, let us define the function $V_L$ as
\begin{equation}\label{V}
V_L(t) := y_x(L(t), t)\quad\  \forall t \in [0, T].
\end{equation}
Obviously, $V_L \in C^\kappa([0, T ])$. Furthermore, from the estimates in the
proofs of Theorems 10.1 and 11.1 in \cite{LSU}, we also have
\begin{equation}\label{e3.1}
\|V_L\|_{C^\kappa([0, T ])}\leq C\|y\|_{C(\overline{R_L})},
\end{equation}
where $C$ is a positive constant depending only on $\alpha$, $N_L$, $L_\ast$, $B$, $\omega$ and $T$.

On the other hand, from the estimates $(i)$ and $(ii)$ in Theorem \ref{t2.1}, we deduce
\begin{equation}\label{i3.2}
\|y\|_{W^{2,1}_2(R_L)}\leq C(\|v\|_{L^2(\omega\times (0,T))}+\|y_0\|_{H^1_\alpha(0,L_0)}),
\end{equation}
where $C$ is a positive constant depending only on $\alpha$, $N_L$, $L_\ast$, $B$, $b'$, $L_0$ and $T$.

Let us now recall the known embedding result
$W^{2,1}_2(R_L)\hookrightarrow C^{1/2,1/4}(\overline{R_L})$, associated with the estimate
\begin{equation}\label{i3.3}
\|y\|_{C^{1/2,1/4}(\overline{R_L})}\leq C\|y\|_{W^{2,1}_2(R_L)},
\end{equation}
where $C$ is a positive constant depending only on $N_L$, $L_\ast$, $B$, $b'$ and $T$.

\eqref{i3.2}, together with \eqref{i3.3}, indicates
\begin{equation}\label{i3.4}
\|y\|_{C(\overline{R_L})}\leq C(\|v\|_{L^2(\omega\times (0,T))}+\|y_0\|_{H^1_\alpha(0,L_0)}).
\end{equation}
Combining  (\ref{e3.1}) with (\ref{i3.4}), we obtain
\begin{equation*}
\|V_L\|_{C^\kappa([0, T ])}\leq C(\|v\|_{L^2(\omega\times (0,T))}+\|y_0\|_{H^1_\alpha(0,L_0)}).
\end{equation*}
Moreover, notice that $v$ satisfies (\ref{cost1}) or (\ref{cost2}). It follows that
\begin{equation}\label{i3.5}
\|V_L\|_{C^\kappa([0, T ])}\leq C\|y_0\|_{H^1_\alpha(0,L_0)},
\end{equation}
where $C$ is a positive constant depending only on $\alpha$, $N_L$, $L_\ast$, $B$, $\omega$ and $T$, and is non-decreasing with respect to $N_L$.

\section{Proof of Theorem \ref{t1.1}}\label{sec3}
To prove Theorem \ref{t1.1}, we will first show that there exists $\varepsilon > 0$ such that, if $\|y_0\|_{H^1_\alpha(0,L_0)}\leq \varepsilon$, then for any $\beta > 0$, there exist uniformly bounded $L_\beta \in C^1([0, T ])$ and control-pairs $(v_\beta , y_\beta )$ satisfying
(\ref{1.1}), (\ref{1.2}) and (\ref{d2}). This can be achieved by applying Schauder's fixed point theorem. For this purpose, let $\beta>0$ and $R>0$ be given. Set
\begin{equation*}
\mathscr M := \{l \in C^1([0, T ]): L_\ast\leq l(t) \leq B,\
l(0) = L_0,\ N_l := \|l'\|_\infty \leq R\}.
\end{equation*}
Clearly, $\mathscr M$ is a non-empty, bounded, closed, and convex subset  of $C^1([0, T ])$. Next, we introduce the mapping $\Lambda_\beta: \mathscr M  \rightarrow C^1([0, T ])$, defined as follows: for each $l\in \mathscr M$, $\Lambda_\beta (l)$ is given by
\begin{equation}\label{iS}
\Lambda_\beta (l)(t)= L_0-\int^t_0l^\alpha(\tau)y_x(l(\tau), \tau) d\tau,
\quad\ t\in[0, T],
\end{equation}
where $y$ is the solution to the linear system
\begin{equation}\label{e3.8}
\left\{\begin{array}{ll}
y_t-(x^\alpha y_x)_x=1_\omega v, & (x,t)\in Q_l,\\[2mm]
\left\{\begin{array}{ll}
y(0,t)=0, & \text{if}\quad \alpha\in[0,1),\\[2mm]
\lim\limits_{x\rightarrow0^+}x^\alpha y_x(x,t)=0, &\text{if}\quad \alpha\in[1,2),\\[2mm]
\end{array}\right.&\\[2mm]
y(l(t),t)=0, &t\in(0,T),\\[2mm]
y(x,0)=y_0(x), &x\in (0,L_0),
\end{array}\right.
\end{equation}
and $v$ is the $\beta$-control constructed via the variational approach presented in Section \ref{s2.5}, such that the solution $y$ satisfies $\|y(x,T)\|_{L^2(0,l(T))}\leq\beta$.

In order to apply Schauder's theorem to $\Lambda_\beta$ in $\mathscr M$, we shall verify that $\Lambda_\beta$ satisfies the following properties:

$(1)$ $\Lambda_\beta$ is continuous. Assume that $l_n\in\mathscr M$ and $l_n\rightarrow l$ in $C^1([0, T ])$ as $n\rightarrow\infty$. We aim to prove that $\Lambda_\beta(l_n)\rightarrow\Lambda_\beta(l)$ (which is equivalent to proving that every subsequence of $\Lambda_\beta(l_n)$ has a further subsequence converging to $\Lambda_\beta(l)$. Let $(v^n,y^n)$ denote the control-state pair of \eqref{e3.8} associated with $l_n$. More precisely, $v^n=\hat{\varphi}^n|_{\omega\times\scriptscriptstyle(0,T)}$, where $\hat{\varphi}^n$ is the solution to (\ref{bf}) for $L = l_n$ and $\varphi_{\scriptscriptstyle T}=\hat{\varphi}^n_{\scriptscriptstyle T}$, and $\varphi^n_{\scriptscriptstyle T}$ is the unique minimizer of $J_\beta(\cdot,l_n)$ (see (\ref{j})). From (\ref{2.26}), we see that for any given $\beta>0$, $\|\hat{\varphi}^n_{\scriptscriptstyle T}\|_{L^2(0,l_n(T))}$ is uniformly bounded with respect to $n$. Thus, there exist a subsequence of $\{\hat{\varphi}^n_{\scriptscriptstyle T}\}_{n=1}^\infty$ (denoted also
by itself for convenience) and a function $\hat{\varphi}_{\scriptscriptstyle T}\in L^2(0,l(T))$ such that their extensions by zero converge weakly in $L^2(0, B)$ to the extension by zero of $\hat{\varphi}_{\scriptscriptstyle T}$.
Moreover, it can be shown that $\hat{\varphi}_{\scriptscriptstyle T}$ is the unique minimizer of $J_\beta(\cdot,l)$, which implies that the sequence $\{\hat{\varphi}^n_{\scriptscriptstyle T}\}_{n=1}^\infty$ converges strongly to $\hat{\varphi}_{\scriptscriptstyle T}$ (see more details in \cite{FLde}).
An immediate consequence is that the controls $v^n$ converge strongly to the control $v=\hat{\varphi}|_{\omega\times\scriptscriptstyle(0,T)}$, where $\hat{\varphi}$ is the solution to (\ref{bf}) for $L = l$ and $\varphi_{\scriptscriptstyle T}=\hat{\varphi}_{\scriptscriptstyle T}$. Also, the states $y^n$ converge strongly to the associated state $y$. Finally, from the estimates in Section \ref{sec3.1}, it follows that $y^n_x(l_n(\cdot),\cdot)\rightarrow y_x(l(\cdot),\cdot)$ strongly in $C([0, T ])$ and $\Lambda_\beta(l_n)\rightarrow\Lambda_\beta(l)$ strongly in $C^1([0,T])$.
Hence, $\Lambda_\beta$ is continuous.

$(2)$ The mapping $\Lambda_\beta$ is compact. In fact, $\Lambda_\beta(l)$ with
$l\in \mathscr M$ belong to a fixed compact set of $C^1([0, T ])$ (independent of
$\beta$). This can be derived from the Arzel\`{a}-Ascoli theorem. We now verify the two key conditions required by the theorem for the family $ \{\Lambda_\beta(l)\}_{l \in M} $, namely uniform boundedness and equicontinuity. First, by the definitions of $\Lambda_\beta(l)$ and $V_l$ given in (\ref{iS}) and (\ref{V}), respectively, we have
\begin{equation*}
\|\Lambda_\beta (l)'\|_{C^\kappa([0, T ])}=\|l^\alpha V_l\|_{C^\kappa([0, T ])}\leq \|l^\alpha\|_{C^\kappa([0, T ])}\| V_l\|_{C^\kappa([0, T ])}.
\end{equation*}
Notice that for any $\alpha>0$, $\kappa\in(0,1)$ and $l\in \mathscr M$,
\begin{eqnarray*}
\|l^\alpha\|_{C^\kappa([0, T ])}\!\!\!\!\!\!\!\!\!&&=\|l^\alpha\|_{C([0, T ])}+\sup\limits_{t_1\neq t_2} \frac{|l^\alpha(t_1)-l^\alpha(t_2)|}{|t_1-t_2|^\kappa}\\[2mm]
&&=\|l^\alpha\|_{C([0, T ])}+\sup\limits_{t_1\neq t_2} \frac{|l^\alpha(t_1)-l^\alpha(t_2)|}{|t_1-t_2|}{|t_1-t_2|^{1-\kappa}}\\[2mm]
&&\leq B^\alpha+\|(l^\alpha)'\|_{C([0, T ])}T^{1-\kappa}\\[2mm]
&&\leq B^\alpha+\alpha\max\{L_\ast^{\alpha-1},B^{\alpha-1}\}\|l'\|_{C([0, T ])}T^{1-\kappa}\\[2mm]
&&\leq B^\alpha+\alpha\max\{L_\ast^{\alpha-1},B^{\alpha-1}\}RT^{1-\kappa}.
\end{eqnarray*}
Combining the two inequalities above with the estimate
(\ref{i3.5}), we conclude that
\begin{eqnarray*}
\|\Lambda_\beta (l)'\|_{C^\kappa([0, T ])}\leq C\|y_0\|_{H^1_\alpha(0,L_0)}.
\end{eqnarray*}
This means the image $\Lambda_\beta(l)$ belongs to a fixed bounded set
of $C^{1+\kappa}([0, T ])$ for some $\kappa \in(0, 1/2]$ that is independent of $l$.

On the other hand, for any $t_1,t_2\in [0,T]$, it holds that
\begin{eqnarray*}
\big|\Lambda_\beta (l)(t_1)-\Lambda_\beta (l)(t_2)\big|\!\!\!\!\!\!\!\!\!&&\leq\|\Lambda_\beta (l)'\|_{C^\kappa([0, T ])}{|t_1-t_2|}
\end{eqnarray*}
and
\begin{eqnarray*}
\big|\Lambda_\beta (l)'(t_1)-\Lambda_\beta (l)'(t_2)\big|\!\!\!\!\!\!\!\!\!&&\leq\|\Lambda_\beta (l)'\|_{C^\kappa([0, T ])}{|t_1-t_2|^{\kappa}}.
\end{eqnarray*}
This means the images $\Lambda_\beta(l)$ and $\Lambda_\beta(l)'$ are equicontinuous.

It follows from the Arzel\`{a}-Ascoli theorem that  $\{\Lambda_\beta(l)\}_{l\in \mathscr M}$ is a relatively compact subset of $C^1([0, T ])$ (i.e., its closure is compact ). Hence, $\Lambda_\beta$ is compact.

$(3)$ $\Lambda_\beta$ maps $\mathscr M$ into itself, provided that $\|y_0\|_{H^1_\alpha(0,L_0)}$ is small enough (uniformly with
respect to $\beta$). This will be done using (\ref{i3.5}). First, by the definition of $\Lambda_\beta(l)$ in (\ref{iS}),
it is clear that $\Lambda_\beta(l) \in C^1([0, T ])$ and $\Lambda_\beta(l)(0) = L_0$. Furthermore,  differentiating both sides of equation (\ref{iS}) with respect to $t$, we get
\begin{equation*}
\big|\Lambda_\beta(l)'(t)\big|=\big|l^\alpha(t)y_x(l(t),t)\big|.
\end{equation*}
Note that $\alpha>0$ and $L_\ast\leq l(t)\leq B$ for all $t\in[0,T]$. It follows that
\begin{equation*}
\big|\Lambda_\beta(l)'(t)\big|\leq\big|l^\alpha(t)\big|\big|y_x(l(t),t)\big|\leq B^\alpha \|V_l\|_{C([0, T ])}.
\end{equation*}
Applying (\ref{i3.5}),  we obtain
\begin{equation*}
\big|\Lambda_\beta(l)'(t)\big|\leq C_{\scriptscriptstyle V}\|y_0\|_{H^1_\alpha(0,L_0)}
\end{equation*}
and moreover from (\ref{iS}),
\begin{equation*}
\big|\Lambda_\beta(l)(t)-L_0\big|\leq C_{\scriptscriptstyle V} T\|y_0\|_{H^1_\alpha(0,L_0)},
\end{equation*}
where $C_{\scriptscriptstyle V}$ denotes a positive constant depending only on $R$, $\alpha$, $N_l$, $L_\ast$, $B$, $\omega$ and $T$. Therefore, if we choose
\begin{equation*}
\|y_0\|_{H^1_\alpha(0,L_0)}\leq  \max\bigg\{\frac{R}{C_{\scriptscriptstyle V}},\frac{B-L_0}{C_{\scriptscriptstyle V}T},\frac{L_0-L_\ast}{C_{\scriptscriptstyle V}T}\bigg\},
\end{equation*}
we will have
\begin{equation*}
 |\Lambda_\beta(l)'(t)| \leq R \quad \text{and} \quad  L_\ast\leq \Lambda_\beta(l)(t) \leq B, \quad \forall t \in[0, T ],
\end{equation*}
whence  the desired property $ \Lambda_\beta(l) \in\mathscr M$ holds.

Now, we can apply Schauder's theorem to $\Lambda_\beta$ in $\mathscr M$ and deduce that there exists a fixed point $L_\beta$ of $\Lambda_\beta$ for each $\beta >0$. Then, it is clear that $L_\beta=\Lambda_\beta(L_\beta)$ satisfies, together with  $v_\beta$ and $y_\beta$, (\ref{1.1}), (\ref{1.2}) and (\ref{d2}).
Moreover, $L_\beta$ and $v_\beta$ are uniformly bounded in $C^{1+\kappa}([0, T ])$ and
$L^2(\omega\times(0,T))$, respectively.
Consequently, at least for a subsequence, one has
\begin{equation*}
L_\beta\rightarrow L \quad \text{strongly in }   \quad   C^1([0, T ])
\end{equation*}
and
\begin{equation*}
\ v_\beta \rightharpoonup v \quad \text{weakly in }   \quad   L^2(\omega\times(0,T))
\end{equation*}
as $\beta\rightarrow0$. Obviously, $L\in\mathscr M$. Also, it is clear that the solution to (\ref{1.1}) satisfies (\ref{1.2}) and (\ref{d1}).

Hence, Theorem 1.1 is proved.
\endpf

\section{Proof of Theorem \ref{t2.2} (Carleman estimate)}\label{carle}

To prove Theorem \ref{t2.2}, we borrow some ideas from \cite{4}. The whole proof is divided into four  parts.

\subsection{Reformulation of the problem}

First, we reformulate (\ref{2.8}). To this aim, let $\sigma(x,t) = \theta(t)\Phi(x)$, where
$$\Phi(x)>0\ \ \forall x \in[0,B] \quad \  \text{and} \quad \  \theta(t)\rightarrow \infty \ \ \text{as}\ \ t \rightarrow 0^+,\ T^-.$$
For any $s > 0$, define
$$z( x,t) = e^{-s\sigma(x,t)}\varphi(x, t),$$
where $\varphi$ is a solution of (\ref{2.8}). For all $n\in \mathbb N$, we observe that
\begin{equation}\label{z}
 \theta^nz = 0\quad \text{and}\quad z_x = 0\quad \text{at time}\ t = 0\ \text{and}\ t = T.
\end{equation}
Moreover, $z$ satisfies
\begin{equation}\label{eqz}
\left\{\begin{array}{ll}
(e^{s\sigma}z)_{t}+(x^\alpha (e^{s\sigma}z)_x)_x=g, & (x,t)\in Q_L,\\[2mm]
\left\{\begin{array}{ll}
z(0,t)=0, & \text{if}\quad \alpha\in[0,1),\\[2mm]
\lim\limits_{x\rightarrow0^+}(x^\alpha z_x)(x, t) = \lim\limits_{x\rightarrow0^+}-s(x^\alpha\sigma_xz)(x, t), &\text{if}\quad \alpha\in[1,2),\\[2mm]
\end{array}\right. &\\[2mm]
z(L(t),t)=0, &t\in(0,T).
\end{array}\right.
\end{equation}
This equation can be rewritten as:
\begin{equation*}
P_sz:=P^+_sz + P^-_sz =ge^{-s\sigma},
\end{equation*}
where
\begin{equation}\label{Ps}
\begin{array}{ll}
&P^+_sz:= (x^\alpha z_x)_x + s^2\sigma_x^2 x^\alpha z+s \sigma_t z,
\\[2mm]
&P^-_sz:= z_t +2s\sigma_xx^\alpha z_x + s(\sigma_xx^\alpha)_xz.
\end{array}
\end{equation}
Furthermore, we have
\begin{equation}\label{Psz}
\|ge^{-s\sigma}\|^2=\|P^+_sz\|^2+ \|P^-_sz\|^2+2\langle P^+_sz, P^-_sz\rangle\geq 2\langle P^+_sz, P^-_sz\rangle,
\end{equation}
where $\|\cdot\|$ and $\langle \cdot, \cdot\rangle$ denote the canonical norm and scalar product in $L^2(Q_L)$.

\subsection{Computation of the scalar product}

We now calculate the scalar product in $L^2(Q_L)$ of $Uz$ and $Vz$ as follows.
\begin{lemma}\label{lid}
Solutions to (\ref{eqz}) satisfy the identity:
\begin{equation*}
\begin{array}{ll}
&\displaystyle\langle P^+_sz, P^-_sz\rangle=Q_1+Q_2+Q_3+Q_4+R+F,
\end{array}
\end{equation*}
where
\begin{equation}\label{elid}
\begin{array}{ll}
&\displaystyle Q_1:= - s\iint_{Q_L}x^{2\alpha-1}(\alpha \sigma_x+2x\sigma_{xx}) z^2_x,\\[4.5mm]
&\displaystyle Q_2:= -s^3\iint_{Q_L}  x^{2\alpha-1}(\alpha \sigma_x+2x\sigma_{xx})\sigma^2_xz^2,\\[4.5mm]
&\displaystyle Q_3:=- \iint_{Q_L}\Big(\frac{1}{2}s \sigma_{tt}+2s^2\sigma_x\sigma_{tx} x^\alpha \Big)z^2,\\[4.5mm]
&\displaystyle Q_4:=-s\iint_{Q_L} x^\alpha (\sigma_xx^\alpha)_{xx}zz_x,\\[4.5mm]
&\displaystyle R:=\int^T_{0} \Big[s\sigma_x(L(t),t)L^{\alpha}(t)-\frac{1}{2}L'(t)\Big]L^\alpha(t)z^2_x(L(t),t),\\[4.5mm]
&\displaystyle  F:=-\int^T_{0}\Big\{ x^\alpha z_xz_t+s^2\sigma_t\sigma_xx^\alpha z^2+s^3\sigma_x^3x^{2\alpha} z^2\\[4.5mm]
&\displaystyle\qquad\qquad\qquad+s\sigma_x(x^\alpha z_x)^2+s x^\alpha (\sigma_xx^\alpha)_xzz_x\Big\}_{|x=0}.
\end{array}
\end{equation}
\end{lemma}
\emph {Proof of Lemma \ref{lid}.} By the definitions of $P^+_sz$ and $P^-_sz$ given in \eqref{Ps}, the scalar product is decomposed as
\begin{equation*}
\langle P^+_sz, P^-_sz\rangle=I_1+I_2+I_3+I_4,
\end{equation*}
where
\begin{eqnarray*}
I_1\!\!\!\!\!\!\!\!\!&&:= \big\langle z_t, (x^\alpha z_x)_x + s^2\sigma_x^2 x^\alpha z+s \sigma_t z\big\rangle,
\\[2mm]
I_2\!\!\!\!\!\!\!\!\!&&:= \big\langle s \sigma_t z, 2s\sigma_xx^\alpha z_x + s(\sigma_xx^\alpha)_xz\big\rangle,
\\[2mm]
I_3\!\!\!\!\!\!\!\!\!&&:= \big\langle s^2\sigma_x^2 x^\alpha z, 2s\sigma_xx^\alpha z_x + s(\sigma_xx^\alpha)_xz\big\rangle,
\\[2mm]
I_4\!\!\!\!\!\!\!\!\!&&:= \big\langle (x^\alpha z_x)_x, 2s\sigma_xx^\alpha z_x + s(\sigma_xx^\alpha)_xz\big\rangle.
\end{eqnarray*}

Next, we compute $I_1$ to $I_4$ sequentially. \emph{First term} $I_1$:
\begin{eqnarray*}
I_1\!\!\!\!\!\!\!\!\!&&= \iint_{Q_L} (x^\alpha z_x)_xz_t + s^2\sigma_x^2 x^\alpha zz_t+s \sigma_t zz_t
\\[2mm]
&&= \iint_{Q_L} (x^\alpha z_xz_t)_x-x^\alpha z_x z_{xt} + (s^2\sigma_x^2 x^\alpha +s \sigma_t)\Big(\text{\small$\frac{z^2}{2}$}\Big)_t
\\[2mm]
&&= \iint_{Q_L} (x^\alpha z_xz_t)_x-\Big(\text{\small$\frac{1}{2}x^\alpha z^2_x$}\Big)_t + (s^2\sigma_x^2 x^\alpha +s \sigma_t)\Big(\text{\small$\frac{z^2}{2}$}\Big)_t.
\end{eqnarray*}
Integrating by parts and applying boundary conditions, we get
\begin{eqnarray*}
I_1\!\!\!\!\!\!\!\!\!&&=\bigg[\int^{L(t)}_0(s^2\sigma_x^2 x^\alpha +s \sigma_t)\frac{z^2}{2}-\frac{1}{2}x^\alpha z^2_x\bigg]^T_0 +\int_{L} x^\alpha z_xz_tn_1\\[2mm]
&&\quad-\int^T_{0}\Big\{ x^\alpha z_xz_t\Big\}_{|x=0}-\int_{L}\frac{1}{2}x^\alpha z^2_x n_2- \iint_{Q_L}\frac{1}{2}(s^2\sigma_x^2 x^\alpha +s \sigma_t)_tz^2,
\end{eqnarray*}
where $n_1$ and $n_2$ denote the spatial and temporal components of the unit exterior normal vector $\mathbf n$ along boundary $L$, respectively.

By (\ref{z}), the terms integrated in time vanish. Hence
\begin{eqnarray}\label{-i1}
I_1= \int_{L} x^\alpha z_xz_tn_1-\int^T_{0}\Big\{ x^\alpha z_xz_t\Big\}_{|x=0}-\int_{L}\frac{1}{2}x^\alpha z^2_x n_2- \iint_{Q_L}\frac{1}{2}(s^2\sigma_x^2 x^\alpha +s \sigma_t)_tz^2.
\end{eqnarray}
Moreover, since the moving boundary of $Q_L$ is characterized by $x=L(t)$, $L\in C^1([0,T])$, we can formulate the unit exterior normal vector as follows
\begin{eqnarray}\label{n}
\mathbf n=(n_1,n_2)=\frac{1}{\sqrt{1+[L'(t)]^2}}(1,-L'(t)),\quad t \in (0, T).
\end{eqnarray}
On the other hand, differentiating $z(L(t), t) = 0$ with respect to $t$ yields $z_t(L(t),t)=-L'(t)z_x(L(t),t)$. Substituting this and (\ref{n}) into \eqref{-i1}, we derive
\begin{eqnarray}\label{i1}
I_1= -\int^T_0\frac{1}{2}L^\alpha(t)L'(t)z^2_x(L(t),t)-\int^T_{0} \Big\{x^\alpha z_xz_t\Big\}_{|x=0}- \iint_{Q_L}\frac{1}{2}(s^2\sigma_x^2 x^\alpha +s \sigma_t)_tz^2.
\end{eqnarray}

\emph{Second term} $I_2$:
\begin{eqnarray*}
I_2\!\!\!\!\!\!\!\!\!&&= s^2\iint_{Q_L} 2\sigma_t\sigma_xx^\alpha zz_x +\sigma_t(\sigma_xx^\alpha)_xz^2= s^2\iint_{Q_L} \sigma_t\sigma_xx^\alpha(z^2)_x+\sigma_t(\sigma_xx^\alpha)_xz^2\notag
\\[2mm]
&&= s^2\iint_{Q_L}(\sigma_t\sigma_xx^\alpha z^2)_x-(\sigma_t\sigma_xx^\alpha)_x z^2+\sigma_t(\sigma_xx^\alpha)_xz^2= s^2\iint_{Q_L}(\sigma_t\sigma_xx^\alpha z^2)_x-\sigma_{tx}\sigma_xx^\alpha z^2.
\end{eqnarray*}
Therefore,
\begin{eqnarray}\label{i2}
I_2=-s^2\int_0^T\Big\{\sigma_t\sigma_xx^\alpha z^2\Big\}_{|x=0}-s^2\iint_{Q_L}\sigma_{tx}\sigma_xx^\alpha z^2.
\end{eqnarray}

\emph{Third term} $I_3$:
\begin{eqnarray*}
I_3\!\!\!\!\!\!\!\!\!&&= s^3  \iint_{Q_L}  2\sigma_x^3x^{2\alpha}z z_x + \sigma_x^2 x^\alpha(\sigma_xx^\alpha)_xz^2= s^3  \iint_{Q_L}  \sigma_x^3x^{2\alpha} (z^2)_x + \sigma_x^2 x^\alpha(\sigma_xx^\alpha)_xz^2\notag
\\[2mm]
&&= s^3  \iint_{Q_L}  (\sigma_x^3x^{2\alpha} z^2)_x -(\sigma_x^3x^{2\alpha})_x z^2+ \sigma_x^2 x^\alpha(\sigma_xx^\alpha)_xz^2=s^3  \iint_{Q_L}  (\sigma_x^3x^{2\alpha} z^2)_x - (\sigma_x^2 x^\alpha)_x\sigma_xx^\alpha z^2.
\end{eqnarray*}
Thus,
\begin{equation}\label{i3}
\begin{array}{ll}
I_3\!\!\!&\displaystyle= -s^3  \int^T_0 \Big\{ \sigma_x^3x^{2\alpha} z^2\Big\}_{|x=0} -s^3  \iint_{Q_L} (\sigma_x^2 x^\alpha)_x \sigma_xx^\alpha z^2\\[4.5mm]
&\displaystyle= -s^3  \int^T_0 \Big\{ \sigma_x^3x^{2\alpha} z^2\Big\}_{|x=0} -s^3  \iint_{Q_L} x^{2\alpha-1}(\alpha \sigma_x+2x\sigma_{xx})\sigma^2_xz^2.
\end{array}
\end{equation}

\emph{Last term} $I_4$:
\begin{eqnarray*}
I_4\!\!\!\!\!\!\!\!\!&&=s\iint_{Q_L} 2\sigma_xx^\alpha z_x(x^\alpha z_x)_x + (x^\alpha z_x)_x(\sigma_xx^\alpha)_xz= s\iint_{Q_L} \sigma_x\big((x^\alpha z_x)^2\big)_x + (x^\alpha z_x)_x(\sigma_xx^\alpha)_xz\\[2mm]
&&= s\iint_{Q_L} \sigma_x\big((x^\alpha z_x)^2\big)_x + \big(x^\alpha z_x(\sigma_xx^\alpha)_xz\big)_x-x^\alpha (\sigma_xx^\alpha)_{xx}zz_x-x^\alpha(\sigma_xx^\alpha)_x z^2_x\\[2mm]
&&= s\iint_{Q_L} \big(\sigma_x(x^\alpha z_x)^2\big)_x+ \big(x^\alpha z_x(\sigma_xx^\alpha)_xz\big)_x -2\sigma_{xx}(x^\alpha z_x)^2-x^\alpha (\sigma_xx^\alpha)_{xx}zz_x-\alpha x^{2\alpha-1} \sigma_xz^2_x.
\end{eqnarray*}
By integration by parts, we get
\begin{eqnarray*}
I_4\!\!\!\!\!\!\!\!\!
&&= s\int_{L} \sigma_x(x^\alpha z_x)^2n_1-s\int^T_0 \Big\{\sigma_x(x^\alpha z_x)^2\Big\}_{|x=0}-s\int^T_0 \Big\{x^\alpha z_x(\sigma_xx^\alpha)_xz\Big\}_{|x=0}\\[2mm]
&&\quad-s\iint_{Q_L}2\sigma_{xx}(x^\alpha z_x)^2+x^\alpha (\sigma_xx^\alpha)_{xx}zz_x+\alpha x^{2\alpha-1} \sigma_xz^2_x,
\end{eqnarray*}
where $n_1$ is given in \eqref{n}.
Consequently,
\begin{equation}\label{i4}
\begin{array}{ll}
I_4\!\!\!&\displaystyle= s\int^T_{0} \sigma_x(L(t),t)L^{2\alpha}(t) z^2_x(L(t),t)-s\int^T_0 \Big\{\sigma_x(x^\alpha z_x)^2+ x^\alpha z_x(\sigma_xx^\alpha)_xz\Big\}_{|{x=0}}\\[4.5mm]
&\displaystyle\quad-s\iint_{Q_L}x^{2\alpha-1}(\alpha \sigma_x+2x\sigma_{xx}) z^2_x -s\iint_{Q_L} x^\alpha (\sigma_xx^\alpha)_{xx}zz_x.
\end{array}
\end{equation}

Finally, Lemma \ref{lid}  follows from (\ref{i1})--(\ref{i4}). \endpf

\subsection{Lower bound estimates for integrals}

To begin with, let us define
\begin{eqnarray}\label{theta}
\theta(t)=\frac{ 1}{\big(t(T-t)\big)^{\scriptscriptstyle4}}   \quad\quad  \forall   t\in(0,T).
\end{eqnarray}

Recall that $\omega=(a,b)$. Let $\omega_{\scriptscriptstyle0}=(c,d)$ with $a<c<d<b$ and consider a cut-off function {\small $\xi\in {\textstyle C}^{\scriptscriptstyle3} (\mathbb R)$} defined as
\begin{eqnarray*}
\xi(x)=
\begin{cases}
1, &x\in(-\infty,c),\\[1mm]
\rho\big((d-x)/(d-c)\big), &x\in[c,d],\\[1mm]
0, &x\in(d,+\infty),
\end{cases}
\end{eqnarray*}
where the polynomial
\begin{equation*}
\rho(\lambda)=\frac{\scriptstyle\raisebox{-0.4ex}{{\footnotesize67}}}{\scriptstyle\raisebox{0.4ex}{\footnotesize3}}{\textstyle\lambda^{\scriptscriptstyle4}}-\frac{\scriptstyle\raisebox{-0.3ex}{\footnotesize103}}{\scriptstyle\raisebox{0.3ex}{\footnotesize2}}\lambda^{\scriptscriptstyle5}+\text{\small$42$}{\lambda}^{\scriptscriptstyle6}-\frac{\scriptstyle\raisebox{-0.3ex}{\footnotesize71}}{\scriptstyle\raisebox{0.3ex}{\footnotesize6}}\lambda^{\scriptscriptstyle7}.
\end{equation*}
It is easy to see that $\xi$ satisfies the derivative bounds
\begin{eqnarray*}
{\textstyle\xi}^{\scriptscriptstyle(k)}(x)\leq\frac{\text{\small $C$}}{ (d-c)^k},\quad\quad\textstyle \text{\small $k = 1, 2, 3,$}
\end{eqnarray*}
for some constant {\small $C>0$} independent of the interval length $d - c$.

Next, define the weight function {\small $\Phi:  [0,{\textstyle B}]\rightarrow\mathbb R_+$} via the convex combination:
\begin{equation}\label{sgm}
\Phi(x)=\xi(x)\phi(x)+\big(1-\xi(x)\big)\psi(x),
\end{equation}
where
\begin{eqnarray*}
\phi(x)=\begin{cases}
\frac{\text{\small${d}^{\scriptstyle2-\alpha}-{ x}^{\scriptstyle 2-\alpha}$}}{{( 2-\alpha)}^2}, &\text{ \small$0\leq\alpha<2,\ \alpha\neq1,$}\\[2mm]
e^d-e^x, &\text{\small$\alpha=1,$}
\end{cases}
\end{eqnarray*}
and
\begin{eqnarray*}
 \psi(x)=e^{2\|\eta\|_{C^{\raisebox{0.05ex}{\scalebox{0.4}{$0$}}}([0,B])}}-e^{\eta(x)},
\end{eqnarray*}
with $\text{\small$\eta(x)$}=- \frac{ 1}{ d} \text{\footnotesize ($x- B$)+1}$. Moreover, we will take
$$\sigma(x,t)=\theta(t)\Phi(x).$$
\begin{remark}
This particular construction of $\Phi(x)$ ensures the smooth transition between $\varphi(x)$ and $\psi(x)$ over the interval $[c, d]$, while maintaining uniform derivative bounds crucial for subsequent energy estimates.
\end{remark}

In the sequel, we use {\small $C$} to denote a positive constant, depending only on $\alpha$, {\small $N_L$}, $c$, $d$, $L_\ast$, $B$ and $T$, which may be different from one place to another; {\small $C_{\scriptscriptstyle0}$, $C_{\scriptscriptstyle1}$}, etc. are other positive (specific) constants.

With the above choice of $\theta$ and $\Phi$, the distributed and boundary terms from Lemma \ref{lid} can be first computed
and then estimated as follows.

\begin{lemma}\label{le2.2.}
For all $\alpha\in[0,2)$, the distributed terms $Q_i$ $(i=1,2,3,4)$ given in (\ref{elid}) satisfy, for $s$ large enough,
\begin{eqnarray}\label{Q}
\iint_{Q_L}s\theta x^\alpha z^2_x +s^3\theta^3 x^{2-\alpha}z^2\leq C\bigg(\sum\limits^{4}_{i=1}Q_i+\int^T_0\int_{\omega_0}s\theta x^{\alpha} z_x^2+s^3 \theta^3x^{2-\alpha}z^2\bigg),
\end{eqnarray}
where $C$ is a positive constant depending only on $\alpha$, $c$, $d$, $B$ and $T$.
\end{lemma}

\noindent \emph {Proof of Lemma \ref{le2.2.}.} The proof proceeds in three steps.

\noindent{\bf Step 1.} First, we estimate lower bounds for $Q_1$. Since $\sigma(x,t)=\theta(t)\Phi(x)$, it follows that
\begin{eqnarray*}
Q_1:=- s\iint_{Q_L}x^{2\alpha-1}(\alpha \sigma_x+2x\sigma_{xx}) z^2_x=-s\iint_{Q_L} \theta x^{2\alpha-1}(\alpha\Phi_x+2x\Phi_{xx}) z_x^2.
\end{eqnarray*}
Next, based on the definition of $\Phi$ in \eqref{sgm}, we partition the interval $(0,L(t))$ into three parts: $(0,c)$, $[c,d]$ and $(d,L(t))$, and analyze each part as follows.

{\noindent $(i)$ Degenerate interval $(0,c)$:} Since $\Phi(x)=\phi(x)$ in $(0,c)$, a direct calculation shows that
$$
\alpha\Phi_x+2x\Phi_{xx}=\alpha\phi_x+2x\phi_{xx}=-x^{1-\alpha}.
$$
This yields the positive definite term
\begin{eqnarray}\label{4.15}
-\int^c_0x^{2\alpha-1}(\alpha\Phi_x+2x\Phi_{xx}) z^2_xdx=\int^c_0 x^{\alpha}z_x^2dx.
\end{eqnarray}

{\noindent $(ii)$ Transition interval $[c,d]$:} The boundedness $x^{\scriptstyle\alpha-1}\left|\alpha\Phi_x+2x\Phi_{xx}\right|\leq C$ implies
\begin{eqnarray}\label{4.16}
-\int^d_c  x^{2\alpha-1}(\alpha\Phi_x+2x\Phi_{xx}) z_x^2dx\geq-C\int_{\omega_0} x^{\alpha}z_x^2dx.
\end{eqnarray}

{\noindent $(iii)$ Moving endpoint interval $(d,L(t))$:} Here, the exponential structure $\psi(x)$ dominates. Thus,
$$\alpha\Phi_x+2x\Phi_{xx}=\alpha\psi_x+2x\psi_{xx}=- \frac{\raisebox{-0.3ex}{$1$}}{\raisebox{0ex}{$d$}} (\frac{\raisebox{-0.3ex}{$2$}}{\raisebox{0ex}{$d$}} x-\alpha)e^{\eta(x)}.$$
Since $\eta(x)>0$ for all $x\in[d,B]$, we have $e^{\eta(x)}>1$. Then, for any $0\leq\alpha<2$ and any $x\in(d,L(t))$, it holds that
$$(\frac{\raisebox{-0.3ex}{$2$}}{\raisebox{0ex}{$d$}} x-\alpha)e^{\eta(x)}>(2-\alpha)>0.$$
On the other hand, notice that $d<L(t)\leq B$ for all $t\in[0, T ]$. It follows that
$$
x^{2\alpha-1}\geq\min\big\{d^{\alpha-1}, B^{\alpha-1}\big\}x^{\alpha}\quad\ \forall x\in (d,L(t)).
$$
Therefore, we establish the coercivity:
\begin{eqnarray}\label{4.17}
-\int^{L(t)}_dx^{2\alpha-1}(\alpha\Phi_x+2x\Phi_{xx}) z^2_xdx\geq C_0\int^{L(t)}_d x^{\alpha}z_x^2dx,
\end{eqnarray}
where $\text{\small$C_{\scriptscriptstyle0}$}=\frac{\small\raisebox{0ex}{$1$}}{\small\raisebox{0ex}{$d$}}{\text{\small$(2-\alpha)$}}\min\{d^{\scriptscriptstyle\alpha-1}, \text{$B$}^{\scriptscriptstyle\alpha-1}\}$.

Combining \eqref{4.15}, \eqref{4.16} and \eqref{4.17}, we obtain the lower bound for $Q_1$:
\begin{equation}\label{Q3}
\begin{array}{ll}
Q_1\!\!\!&\displaystyle\geq s\int^T_0\int^c_0\theta x^{\alpha}z_x^2-Cs\int^T_0\int_{\omega_0}\theta x^{\alpha} z_x^2+C_0s\int^T_0\int^{L(t)}_d \theta x^{\alpha}z_x^2\\[4.5mm]
&\displaystyle\geq C_1s\iint_{Q_L}\theta x^{\alpha}z_x^2-Cs\int^T_0\int_{\omega_0}\theta x^{\alpha} z_x^2,
\end{array}
\end{equation}
where $C_{\scriptscriptstyle1}=\min\{C_{\scriptscriptstyle0},1\}$.

Similarly, taking into account that $\Phi_x(x)=\phi_x(x)=-\frac{x^{\scriptscriptstyle1-\alpha}}{\scriptscriptstyle{2-\alpha}}$ in $(0,c)$, and $\Phi_x(x)=\psi_x(x)=\frac{1}{d} {\scriptstyle e}^{\scriptscriptstyle\eta(x)}>\frac{1}{d}$ in $(d,L(t))$, we derive
\begin{equation}\label{Q2}
\begin{array}{ll}
Q_2\!\!\!&\displaystyle:=-s^3  \iint_{Q_L}  x^{2\alpha-1}(\alpha \sigma_x+2x\sigma_{xx})\sigma^2_xz^2\\[4.5mm]
&\displaystyle\ =-s^3\iint_{Q_L}\theta^3 x^{2\alpha-1} (\alpha\Phi_x+2x\Phi_{xx})\Phi^2_x z^2\\[4.5mm]
&\displaystyle\ \geq C_2s^3\iint_{Q_L}\theta^3x^{2-\alpha} z^2-Cs^3\int^T_0\int_{\omega_0} \theta^3x^{2-\alpha}z^2,
\end{array}\end{equation}
where $C_{\scriptscriptstyle2}$ is a positive constant depending only on $\alpha$, $d$ and $B$.

\noindent{\bf Step 2.} Let us analyze the third term  $Q_{\scriptscriptstyle3}$:
\begin{eqnarray*}
Q_3:=- \iint_{Q_L}\Big(\frac{1}{2}s \sigma_{tt}+2s^2\sigma_x\sigma_{tx} x^\alpha \Big)z^2=- \iint_{Q_L}\Big(\frac{1}{2}s \theta_{tt}\Phi+2s^2\theta\theta_t\Phi^2_x x^\alpha \Big)z^2.
\end{eqnarray*}
From \eqref{theta}, we observe  that
$$
 |\theta_{tt}|\leq C(T)\theta^{3/2}\quad\quad \text{and}\quad\quad |\theta\theta_t|\leq C(T)\theta^{9/4}\leq C(T)\theta^3.
$$
On the other hand, by the definition of $\Phi$ in \eqref{sgm}, it is easy to check that {\small $|\Phi(x)|\leq C$} and {\small$|\Phi_x(x)|\leq \text{\small$C$}x^{\scriptstyle1-\alpha}$} for all {\small $x\in [0,B]$}. Then, we obtain
\begin{eqnarray}\label{e2.24}
Q_3\leq C\iint_{Q_L}\big(s \theta^{3/2}+s^2\theta^3x^{2-\alpha} \big)z^2.
\end{eqnarray}
It remains to bound the first term on the right-hand side of the above inequality. First, by Theorem \ref{t2.1}, the solution $\varphi$ of \eqref{2.8} belongs to $L^2(0,T;H_\alpha^1(0,L(t)))$. Second, since
$z=e^{-s\sigma}\varphi$, a straightforward calculation shows that $z$ also belongs to $L^2(0,T;H_\alpha^1(0,L(t)))$. Next, for any $\varepsilon> 0$, we write
\begin{equation}\label{qt1}
\begin{array}{ll}
\displaystyle\iint_{Q_L}\theta^{3/2}z^2\!\!\!&\displaystyle=\iint_{Q_L}\big(\theta x^{(\alpha-2)/3}z^2\big)^{3/4}(\theta^{3}x^{2-\alpha}z^2)^{1/4}\\[4.5mm]
&\displaystyle\leq \frac{3\varepsilon}{4}\iint_{Q_L}\theta x^{(\alpha-2)/3}z^2+\frac{1}{4\varepsilon^{3}}\iint_{Q_L}\theta^{3}x^{2-\alpha}z^2.
\end{array}
\end{equation}
At this point, we distinguish the case $\alpha = 1$ from the others, as Hardy's inequality (\ref{ht}) does not hold for $\alpha^\star = 1$.

In the case \text{\small$\alpha\neq1$}, we observe that $x^{(\alpha-2)/3}\leq \text{\small$B$}^{[2(2-\alpha)]/3}x^{\alpha-2}$ (since \text{\small$\alpha<2$}) and apply
Lemma \ref{lht} with $\alpha^\star=\alpha\neq1$ ($z$ satisfies the assumptions of Lemma \ref{lht} for almost
every $t\in(0,T)$ because it belongs to $L^2(0,T;H_\alpha^1(0,L(t)))$) to obtain
\begin{eqnarray}\label{c1}
\iint_{Q_L}\theta x^{(\alpha-2)/3} z^2\leq B^{2(2-\alpha)/3}\iint_{Q_L}\theta x^{\alpha-2}z^2\leq \frac{4B^{[2(2-\alpha)]/3}}{(\alpha-1)^2}\iint_{Q_L}\theta x^{\alpha}z_x^2.
\end{eqnarray}
In the case of $\alpha=1$, we apply Lemma \ref{lht} with $\alpha^\star= 5/3$ and then use the fact that
$x^{5/3} \leq \text{\small$B$}^{2/3}x$ to arrive at a similar conclusion:
\begin{eqnarray}\label{c2}
\begin{array}{ll}
\displaystyle\iint_{Q_L}\theta x^{(\alpha-2)/3} z^2\!\!\!&\displaystyle=\iint_{Q_L}\theta x^{-1/3} z^2\leq \frac{4}{(\alpha^\star-1)^2}\iint_{Q_L}\theta x^{5/3} z_x^2\\[4.5mm]
&\displaystyle\leq 9B^{2/3}\iint_{Q_L}\theta x z_x^2=9B^{2/3}\iint_{Q_L}\theta x^{\alpha} z_x^2.
\end{array}
\end{eqnarray}
In both cases, combining (\ref{qt1}) with (\ref{c1}) or (\ref{c2}), we deduce
\begin{eqnarray}\label{e2.28}
\iint_{Q_L}\theta^{3/2}z^2\!\!\!\!\!\!\!\!\!&&\leq C\varepsilon\iint_{Q_L}\theta x^{\alpha} z_x^2+\frac{1}{4\varepsilon^{3}}\iint_{Q_L}\theta^{3}x^{2-\alpha}z^2.
\end{eqnarray}
Summing up, we obtain by (\ref{e2.24}) and (\ref{e2.28})
\begin{eqnarray}\label{Q1}
Q_3\leq Cs\varepsilon\iint_{Q_L}\theta x^{\alpha} z_x^2+Cs\varepsilon^{-3}\iint_{Q_L}\theta^{3}x^{2-\alpha}z^2+Cs^2\iint_{Q_L}\theta^{3}x^{2-\alpha}z^2.
\end{eqnarray}

\noindent{\bf Step 3.} Next, we estimate the last term $\textstyle Q_{\scriptscriptstyle4}$. Since $\Phi_{\scriptstyle x}x^{\scriptstyle\alpha}=-\frac{\text{\footnotesize$x$}}{\scriptstyle2-\alpha}$ in $(0,c)$, it follows that $(\Phi_{\scriptstyle x}x^{\scriptstyle\alpha})_{\scriptstyle xx}=0$ in $(0,c)$. Hence,
\begin{eqnarray*}
Q_4\!\!\!\!\!\!\!\!\!&&:=- s\iint_{Q_L} x^\alpha (\sigma_xx^\alpha)_{xx}zz_x=-s\iint_{Q_L}\theta x^\alpha (\Phi_xx^\alpha)_{xx}zz_x=-s\int^T_0\!\!\!\int^{L(t)}_c\theta x^\alpha (\Phi_xx^\alpha)_{xx}zz_x.
\end{eqnarray*}
By {\small(\ref{theta})}, we have {\small$|\theta|\leq C(T)\theta^{\scriptstyle3}$}. Additionally, notice that {\small $|(\Phi_{\scriptstyle x}x^{\scriptstyle\alpha})_{\scriptstyle xx}|\leq C$} for  $x\in(c,L(t))$ and that {\small $x^{\scriptstyle\alpha}\leq c_1x$} holds for all $\text{\small$x$}\in[c,\text{\small$B$}]$ with {\small $c_1=\max\{c^{\scriptstyle\alpha-1},\text{\small$B$}^{\scriptstyle\alpha-1}\}$}. We consequently conclude that for any $\varepsilon> 0$,
\begin{eqnarray}\label{Q4}
Q_4\!\!\!\!\!\!\!\!\!&&\leq Cs\int^T_0\!\!\!\int^{L(t)}_c\theta x |zz_x|\leq Cs\Big(\varepsilon\iint_{Q_L}\theta x^\alpha z^2_x +\varepsilon^{-1}\iint_{Q_L}\theta^3 x^{2-\alpha}z^2\Big).
\end{eqnarray}

In summary, for $\varepsilon$ sufficiently small and $s$ sufficiently large, by combining (\ref{Q3}), (\ref{Q2}), (\ref{Q1}) and (\ref{Q4}), we obtain the desired estimate \eqref{Q}.\endpf

\begin{lemma}\label{leRF}
For all $\alpha\in[0,2)$, the boundary
terms $R$ and $F$ given in (\ref{elid}) satisfy, for $s$ large enough,
\begin{eqnarray}\label{RF}
s\int^T_0L^\alpha(t) z^2_x(L(t),t)\leq C\big(R+F\big),
\end{eqnarray}
where $C$ is a positive constant depending only on $\alpha$, $N_L$, $d$, $L_\ast$ and $T$.
\end{lemma}

\emph{Proof of Lemma \ref{leRF}.} Let us first analyze the degenerate boundary terms:
\begin{eqnarray*}
F\!\!\!\!\!\!\!\!\!&&:= -\int^T_{0}\Big\{ x^\alpha z_xz_t+s^2\sigma_t\sigma_xx^\alpha z^2+s^3 \sigma_x^3x^{2\alpha} z^2+s\sigma_x(x^\alpha z_x)^2+s x^\alpha (\sigma_xx^\alpha)_xzz_x\Big\}_{|x=0}.
\end{eqnarray*}

In the case $0\leq\alpha<1$, applying the boundary condition $z(0,t)=0$ (which implies $z_t(0,t)=0$), we obtain
\begin{eqnarray*}
F_{\ 0\leq\alpha<1}= -\int^T_{0}\Big\{s\sigma_x(x^\alpha z_x)^2\Big\}_{|x=0}.
\end{eqnarray*}
Recalling that $\sigma(x,t)=\theta(t)\Phi(x)$ and using \eqref{sgm}, we conclude that
\begin{eqnarray*}\label{bt1}
F_{\ 0\leq\alpha<1}= \frac{1}{2-\alpha}\int^T_{0}\Big\{s\theta x^{1+\alpha} z_x^2\Big\}_{|x=0}.
\end{eqnarray*}

In the case $1\leq\alpha<2$, since $z$ satisfies the boundary condition $\lim\limits_{x\rightarrow0^+}(x^\alpha z_x)(x, t) = \lim\limits_{x\rightarrow0^+}-s(x^\alpha\sigma_xz)(x, t)$, the boundary term $F$ becomes
\begin{eqnarray*}
F_{\ 1\leq\alpha<2}\!\!\!\!\!\!\!\!\!&&=\int^T_{0}\bigg\{sx^\alpha\sigma_x\Big(\text{\small $\frac{z^2}{2}$}\Big)_{\scriptscriptstyle t} -s^2\sigma_t\sigma_xx^\alpha z^2 -2s^3 \sigma_x^3x^{2\alpha} z^2 + s^2 (\sigma_xx^\alpha)_xx^\alpha\sigma_xz^2\bigg\}_{|x=0}\\[2mm]
&&=\int^T_{0}\Big\{-\frac{s}{2}x^\alpha\sigma_{xt}z^2 -s^2\sigma_t\sigma_xx^\alpha z^2 -2s^3 \sigma_x^3x^{2\alpha} z^2 + s^2 (\sigma_xx^\alpha)_xx^\alpha\sigma_xz^2\Big\}_{|x=0}.
\end{eqnarray*}
Using the expression for $\sigma$ and \eqref{sgm}, we get
\begin{eqnarray*}
F_{\ 1\leq\alpha<2}=\int^T_{0}\bigg\{\bigg(\frac{s\theta_t}{2(2-\alpha)}+\frac{2s^2d^{2-\alpha}\theta_t\theta}{(2-\alpha)^3}+ \frac{2s^3\theta^3}{(2-\alpha)^3}x^{2-\alpha}+\frac{s^2\theta^2}{(2-\alpha)^2}\bigg)xz^2\bigg\}_{|x=0}.
\end{eqnarray*}
Since $z\in H^1_\alpha(0,L(t))$ for almost every $t\in(0,T)$, it follows that $xz^2(x,t)\rightarrow0$ as $x\rightarrow0^+$(\cite[Lemma 3.5]{4}). Hence,
\begin{eqnarray*}\label{bt2}
F_{\ 1\leq\alpha<2}=0.
\end{eqnarray*}
In both cases, we obtain
\begin{eqnarray}\label{bt3}
F\geq0.
\end{eqnarray}

We now turn to the moving boundary terms:
\begin{eqnarray*}
R:=\int^T_{0} \Big[s\sigma_x(L(t),t)L^{\alpha}(t)-\frac{1}{2}L'(t)\Big]L^\alpha(t)z^2_x(L(t),t).
\end{eqnarray*}
Substituting $\sigma(x,t)=\theta(t)\Phi(x)$ and \eqref{sgm} into the above formula yields
\begin{eqnarray*}
R\!\!\!\!\!\!\!\!\!&&=\int^T_{0} \Big[s\theta\psi_x(L(t))L^{\alpha}(t)-\frac{1}{2}L'(t)\Big]L^\alpha(t)z^2_x(L(t),t)\notag\\[2mm]
&&=\int^T_{0} \Big[s\theta\frac{1}{d}e^{\eta(L(t))}L^{\alpha}(t)-\frac{1}{2}L'(t)\Big]L^\alpha(t)z^2_x(L(t),t).
\end{eqnarray*}
Notice that $\eta(L(t))\geq1$, $L(t)\geq L_\ast$ and $\theta(t)\geq\theta(T/2)$ for all $t\in[0,T]$. It follows that
\begin{eqnarray*}
R\geq\int^T_{0} \Big(s\theta(T/2)\frac{e}{d}L_\ast^{\alpha}-\frac{1}{2}N_L\Big)L^\alpha(t)z^2_x(L(t),t).
\end{eqnarray*}
Thus, for $s$ large enough,
\begin{eqnarray}\label{emb}
R\geq C s\int^T_{0} L^\alpha(t)z^2_x(L(t),t),
\end{eqnarray}
where $C>0$ depending only on $\alpha$, $N_L$, $d$, $L_\ast$ and $T$.

Combining (\ref{bt3}) with (\ref{emb}), we obtain the desired inequality (\ref{RF}). \endpf

\subsection{Conclusion}
From Lemmas \ref{lid}-\ref{leRF} and \eqref{Psz}, we conclude that for all $0\leq\alpha<2$,
\begin{eqnarray*}
&&\quad\iint_{Q_L}\big(s\theta x^\alpha z^2_x+s^3\theta^3x^{2-\alpha}z^2\big)+s\int^T_0L^\alpha(t) z^2_x(L(t),t)\\[2mm]
&&\leq C\(\int^T_0\int_{\omega_0}\big(s\theta x^\alpha z^2_x+s^3\theta^3x^{2-\alpha}z^2\big)+\iint_{Q_L}g^2e^{-2s\sigma}\).
\end{eqnarray*}
Recalling that $\varphi=e^{s\sigma}z$, we have $\varphi_x=e^{s\sigma}z_x+s\sigma_xe^{s\sigma}z$. This, together with the definition of $\sigma$, yields
\begin{eqnarray*}
s\theta x^\alpha \varphi^2_x+s^3\theta^3x^{2-\alpha}\varphi^2\!\!\!\!\!\!\!\!\!&&\leq s\theta x^\alpha \big(2e^{2s\sigma}z^2_x+2s^2\sigma_x^2e^{2s\sigma}z^2\big)+s^3\theta^3x^{2-\alpha}e^{2s\sigma}z^2\\[2mm]
&&\leq C\(s\theta x^\alpha e^{2s\sigma}z^2_x+s^3\theta^3x^{2-\alpha}e^{2s\sigma}z^2\).
\end{eqnarray*}
Conversely, since $z=e^{-s\sigma}\varphi$, we get $z_x=e^{-s\sigma}\varphi_x-s\sigma_xe^{-s\sigma}\varphi$. Moreover, notice that $\varphi(L(t),t)=0$, which implies $z_x(L(t),t)=e^{-s\sigma(L(t),t)}\varphi_x(L(t),t)$. Finally, we obtain
\begin{eqnarray*}
\iint_{Q_L}\big(s\theta x^\alpha \varphi^2_x+s^3\theta^3x^{2-\alpha}\varphi^2\big)e^{-2s\sigma}+s\int^T_0L^\alpha(t) e^{-2s\sigma(L(t),t)}\varphi^2_x(L(t),t)\\[2mm]
\leq C\(\int^T_0\int_{\omega_0}\big(s\theta x^\alpha \varphi^2_x+s^3\theta^3x^{2-\alpha}\varphi^2\big)e^{-2s\sigma}+\iint_{Q_L}g^2e^{-2s\sigma}\).
\end{eqnarray*}
This completes the proof of Theorem \ref{t2.2}.\endpf

\section{Appendix}\label{sec6}

\subsection{Proof of Lemma \ref{lht} (Hardy's inequalities)}

\emph{First case.} $0 \leq \alpha^\star < 1$. Since $z$ is absolutely continuous on $(0, L(t))$, we have
\begin{eqnarray*}
|z(x)-z(\varepsilon)|^2\!\!\!\!\!\!\!\!\! &&=\bigg(\int^x_\varepsilon z_x(\mu)\mu^{(3-\gamma)/4}\mu^{(-3+\gamma)/4} d\mu\bigg)^2\\[2mm]
&&\leq \bigg(\int^x_\varepsilon z^2_x(\mu)\mu^{(3-\gamma)/2} d\mu\bigg)\bigg(\int^x_\varepsilon \mu^{(-3+\gamma)/2} d\mu\bigg),
\end{eqnarray*}
where we denote $\gamma := 2-\alpha^\star \in(1, 2]$. Letting $\varepsilon\rightarrow 0^+$, we get
\begin{eqnarray*}
|z(x)|^2\leq \bigg(\int^x_0 z^2_x(\mu)\mu^{(3-\gamma)/2} d\mu\bigg)\bigg(\int^x_0 \mu^{(-3+\gamma)/2} d\mu\bigg).
\end{eqnarray*}
Therefore,
\begin{eqnarray*}
&&\int^{L(t)}_0x^{\alpha^\star-2}|z(x)|^2dx\leq \int^{L(t)}_0x^{-\gamma}\bigg(\int^x_0 z^2_x(\mu)\mu^{(3-\gamma)/2} d\mu\bigg)\bigg(\int^x_0 \mu^{(-3+\gamma)/2} d\mu\bigg)dx\\[2mm]
&&\quad\quad\quad\quad\quad\quad\quad\quad\quad\ =\int^{L(t)}_0x^{-\gamma}\bigg(\int^x_0 z^2_x(\mu)\mu^{(3-\gamma)/2} d\mu\bigg)\frac{x^{(\gamma-1)/2}}{(\gamma-1)/2}dx\\[2mm]
&&\quad\quad\quad\quad\quad\quad\quad\quad\quad\ =\frac{2}{\gamma-1}  \int^{L(t)}_0z^2_x(\mu)\mu^{(3-\gamma)/2}\bigg(\int^{L(t)}_\mu x^{(-\gamma-1)/2}dx\bigg) d\mu\\[2mm]
&&\quad\quad\quad\quad\leq\frac{2}{\gamma-1}  \int^{L(t)}_0z^2_x(\mu)\mu^{(3-\gamma)/2}\frac{\mu^{(1-\gamma)/2}}{(\gamma-1)/2} d\mu= \frac{4}{(1-\alpha^\star)^2}  \int^{L(t)}_0\mu^{\alpha^\star}z^2_x(\mu) d\mu .
\end{eqnarray*}

\emph{Second case.} $1 < \alpha^\star < 2$. Denoting $\gamma := 2-\alpha^\star \in(0, 1)$, we have
\begin{eqnarray*}
&&\int^{L(t)}_0x^{\alpha^\star-2}|z(x)|^2dx\leq \int^{L(t)}_0x^{-\gamma}\bigg(\int^{L(t)}_x z^2_x(\mu)\mu^{(3-\gamma)/2} d\mu\bigg)\bigg(\int^{L(t)}_x \mu^{(-3+\gamma)/2} d\mu\bigg)dx\\[2mm]
&&\quad\quad\quad\quad\quad\quad\quad\quad\quad\ \leq\int^{L(t)}_0x^{-\gamma}\bigg(\int^{L(t)}_x z^2_x(\mu)\mu^{(3-\gamma)/2} d\mu\bigg)\frac{x^{(\gamma-1)/2}}{(1-\gamma)/2}dx\\[2mm]
&&\quad\quad\quad\quad\quad\quad\quad\quad\quad\ =\frac{2}{1-\gamma}  \int^{L(t)}_0z^2_x(\mu)\mu^{(3-\gamma)/2}\bigg(\int^\mu_0 x^{(-\gamma-1)/2}dx\bigg) d\mu\\[2mm]
&&\quad\quad\quad\quad\quad \leq\frac{2}{1-\gamma}  \int^{L(t)}_0z^2_x(\mu)\mu^{(3-\gamma)/2}\frac{\mu^{(1-\gamma)/2}}{(1-\gamma)/2} d\mu= \frac{4}{(\alpha^\star-1)^2}  \int^{L(t)}_0\mu^{\alpha^\star}z^2_x(\mu) d\mu .
\end{eqnarray*}
Hence, we complete the proof of Lemma \ref{lht}.\endpf

\subsection{Proof of Lemma \ref{cacc} (Caccioppoli's inequality)}

Consider a smooth function $\zeta: \mathbb R  \rightarrow \mathbb R $ such that
\begin{equation*}
\left\{\begin{array}{ll}
0\leq\zeta(x)\leq1 &\text{for}\ x\in \mathbb R,\\[2mm]
\zeta(x)=1 &\text{for}\ x\in\omega_0,\\[2mm]
\zeta(x)=0 &\text{for}\ x\notin \omega.
\end{array}\right.
\end{equation*}
Then, for all $s >0$,
\begin{eqnarray*}
0\!\!\!\!\!\!\!\!\!&&=\int_0^T\frac{d}{dt}\int_0^{L(t)}\zeta^2 e^{-2s\sigma}\varphi^2=\iint_{Q_L}-2\zeta^2s\sigma_t e^{-2s\sigma}\varphi^2+2\zeta^2e^{-2s\sigma}\varphi\varphi_t\\[2mm]
&&=-2\iint_{Q_L}\zeta^2s\sigma_te^{-2s\sigma}\varphi^2-2\iint_{Q_L}\zeta^2e^{-2s\sigma}\varphi(x^\alpha\varphi_x)_x\\[2mm]
&&=-2\iint_{Q_L}\zeta^2s\sigma_te^{-2s\sigma}\varphi^2+2\iint_{Q_L}(\zeta^2e^{-2s\sigma}\varphi)_x x^\alpha\varphi_x\\[2mm]
&&=-2\iint_{Q_L}\zeta^2s\sigma_te^{-2s\sigma}\varphi^2+2\iint_{Q_L}x^\alpha(\zeta^2e^{-2s\sigma})_x \varphi\varphi_x+\zeta^2e^{-2s\sigma} x^\alpha\varphi^2_x.
\end{eqnarray*}
Hence,
\begin{eqnarray*}
&&2\iint_{Q_L}\zeta^2e^{-2s\sigma} x^\alpha\varphi^2_x=\iint_{Q_L}2\zeta^2s\sigma_te^{-2s\sigma}\varphi^2-2\iint_{Q_L}x^\alpha(\zeta^2e^{-2s\sigma})_x \varphi\varphi_x\\[2mm]
=\!\!\!\!\!\!\!\!\!&&2\iint_{Q_L}\zeta^2s\sigma_te^{-2s\sigma}\varphi^2-2\iint_{Q_L}\big(x^{\alpha/2}\zeta e^{-s\sigma}\varphi_x\big)\bigg(x^{\alpha/2}\frac{(\zeta^2e^{-2s\sigma})_x}{\zeta e^{-s\sigma}}\varphi\bigg)\\[2mm]
\leq\!\!\!\!\!\!\!\!\!&&2\iint_{Q_L}\zeta^2s\sigma_te^{-2s\sigma}\varphi^2+\iint_{Q_L}\big(x^{\alpha/2}\zeta e^{-s\sigma}\varphi_x\big)^2+\iint_{Q_L}\bigg(x^{\alpha/2}\frac{(\zeta^2e^{-2s\sigma})_x}{\zeta e^{-s\sigma}}\varphi\bigg)^2\\[2mm]
=\!\!\!\!\!\!\!\!\!&&2\iint_{Q_L}\zeta^2s\sigma_te^{-2s\sigma}\varphi^2+\iint_{Q_L}\bigg(x^{\alpha/2}\frac{(\zeta^2e^{-2s\sigma})_x}{\zeta e^{-s\sigma}}\varphi\bigg)^2+\iint_{Q_L}\zeta^2 e^{-2s\sigma}x^{\alpha}\varphi^2_x.
\end{eqnarray*}
Therefore,
\begin{eqnarray*}
\iint_{Q_L}\zeta^2 e^{-2s\sigma}x^{\alpha}\varphi^2_x\leq2\iint_{Q_L}\zeta^2s\sigma_te^{-2s\sigma}\varphi^2+\iint_{Q_L}\bigg(x^{\alpha/2}\frac{(\zeta^2e^{-2s\sigma})_x}{\zeta e^{-s\sigma}}\varphi\bigg)^2\leq C(s,T)\int^T_0\int_{\omega}\varphi^2dxdt,
\end{eqnarray*}
which completes the proof of Lemma \ref{cacc}.\endpf

\end{document}